\documentclass[11pt,reqno]{amsart}
\usepackage[a4paper,margin=25mm]{geometry}
\usepackage{amsmath,amssymb,amsthm,mathtools,booktabs,array,longtable,xcolor}
\usepackage[colorlinks=true,linkcolor=blue!45!black,citecolor=blue!45!black,urlcolor=blue!45!black]{hyperref}
\usepackage{microtype}
\hypersetup{pdftitle={On the Haagerup vertex operator algebras},pdfauthor={Chongying Dong, Xiangyu Jiao, Feng Xu}}
\newtheorem{theorem}{Theorem}[section]
\newtheorem{proposition}[theorem]{Proposition}
\newtheorem{lemma}[theorem]{Lemma}
\newtheorem{corollary}[theorem]{Corollary}
\theoremstyle{definition}

\theoremstyle{remark}\newtheorem{remark}[theorem]{Remark}
\newcommand{\C}{\mathbb C}\newcommand{\Z}{\mathbb Z}\newcommand{\Q}{\mathbb Q}
\newcommand{\one}{\mathbf1}\newcommand{\slc}{\mathfrak{sl}_2}
\newcommand{\Rep}{\operatorname{Rep}}
\newcommand{\FPdim}{\operatorname{FPdim}}
\newcommand{\wt}{\operatorname{wt}}\newcommand{\eps}{\epsilon}
\newcommand{\Tc}{\mathcal T}\newcommand{\D}{\mathcal D}
\newcommand{\nn}[1]{\mathopen{:}#1\mathclose{:}}
\numberwithin{equation}{section}
\title[On the Haagerup vertex operator algebras]{On the Haagerup vertex operator algebras}
\author{Chongying Dong}
\address{Department of Mathematics, University of California, Santa Cruz, California, USA}
\email{dong@ucsc.edu}
\author{Xiangyu Jiao}
\address{School of Mathematical Sciences, East China Normal University, Shanghai, China}
\email{xyjiao@math.ecnu.edu.cn}
\author{Feng Xu}
\address{Department of Mathematics, University of California, Riverside, California, USA}
\email{xufeng@math.ucr.edu }
\date{28 October 2026}
\begin{document}
\begin{abstract}
Three constructions of vertex operator algebras realizing the Haagerup center were announced independently and almost simultaneously. The relations between these constructions are studied in this paper. The generators are expressed in a common lattice basis for $E_6\oplus A_2$. It is proved that the constructions of Riesen and Huang give isomorphic vertex operator algebras, and that this algebra and the algebra constructed by Jiao have isomorphic fixed point subalgebras under automorphisms of order two. The isomorphisms are given explicitly on generators. The common fixed point subalgebra contains the parafermion vertex operator algebras $K(\slc,3)$ and $K(\slc,13)$. A general rationality criterion for unitary vertex operator algebras is given in terms of positive linear combinations of conformal vectors and is applied to these constructions. The two vertex operator algebras are simple current extensions of the common fixed point subalgebra and have equivalent module categories as ribbon categories. {The equivalence follows from decomposition with respect to a nondegenerate pointed subcategory with four simple objects in the module category of the common fixed point subalgebra.} It is also proved that the conformal nets associated with these unitary vertex operator algebras are completely rational.
\end{abstract}
\maketitle

\section{Introduction}
The purpose of this paper is to determine the relations between three vertex operator algebra realizations of the Haagerup center. In this introduction, we explain the motivation for this comparison, state the main results, and describe the general rationality criterion used in their proof.

The realization of modular tensor categories by rational vertex operator algebras is an important problem in the representation theory of vertex operator algebras. The Haagerup fusion categories are of particular interest in this problem. They arise in subfactor theory and lie outside the familiar constructions from finite groups and affine Lie algebras. The six-object Haagerup category has a noncommutative fusion ring and hence cannot itself be the module category of a rational vertex operator algebra. Its Drinfeld center is a modular tensor category. The problem is therefore to realize this center as the module category of a rational vertex operator algebra. We denote the six-object Haagerup category by $\mathcal H_3$ and its Drinfeld center by $Z(\mathcal H_3)$.

Evans and Gannon studied the modular data of the Haagerup center and obtained possible character vectors at central charge eight \cite{EG}. These results suggested looking for a realization in the holomorphic lattice vertex operator algebra $V_{E_8}$. The lattice inclusion $E_6\oplus A_2\subset E_8$ gives a conformal inclusion
\[
 B:=V_{E_6\oplus A_2}\subset E:=V_{E_8}.
\]
It provides a setting in which the generators and their products can be described explicitly.

Three independent constructions were announced almost simultaneously in October 2026 by Riesen \cite{Riesen,RiesenLong}, Huang \cite{Huang} and Jiao \cite{Jiao}. All three are given by conformal subalgebras of $B$. The constructions of Riesen and Huang contain commuting affine vertex operator algebras at levels $1$ and $39$, while Jiao's construction contains commuting affine vertex operator algebras at levels $3$ and $13$. We write $L_{\mathfrak g}(k,0)$ for the simple affine vertex operator algebra associated to a simple Lie algebra $\mathfrak g$ at level $k$, and put $\slc=\mathfrak{sl}_2(\C)$.

The character and weight-one calculations in the original constructions give a further reason to compare them \cite{Huang,Jiao}. Their ordinary vacuum characters are equal:
\begin{equation}\label{eq:commoncharacter}
 \operatorname{ch}(q)=q^{-1/3}
 \bigl(1+6q+120q^2+956q^3+6010q^4+\cdots\bigr).
\end{equation}
Here $q$ is the character variable, and the character of a vertex operator algebra $A$ of central charge $c_A$ is $\operatorname{ch}_A(q)=\operatorname{Tr}_A q^{L(0)-c_A/24}$. This is the vacuum component of the $\gamma=0$ Evans--Gannon character vector, where $\gamma$ is the parameter in their character classification \cite{EG}. Their weight-one Lie algebras are all isomorphic to $\slc\oplus\slc$. The corresponding affine vertex operator algebras are, however,
\[
 L_{\slc}(1,0)\otimes L_{\slc}(39,0)
 \quad\text{and}\quad
 L_{\slc}(3,0)\otimes L_{\slc}(13,0),
\]
respectively. Thus the affine levels distinguish two realizations with the same character and the same abstract weight-one Lie algebra.

The main purpose of this paper is to determine the relations between these three constructions. We express all generators in one ordered lattice basis and fix the cocycle. The two additional Virasoro vectors in Riesen's presentation are recovered from Huang's generators. This proves that the Riesen and Huang constructions are isomorphic. We then show that this vertex operator algebra and Jiao's vertex operator algebra have isomorphic fixed point subalgebras under automorphisms of order two.

Let $U_{\rm Riesen}$ denote Riesen's vertex operator algebra, and let $U$ and $V$ denote the algebras defined by Huang and Jiao, respectively. Their generators are given in Section~\ref{sec:lattice}. For a rational, $C_2$-cofinite vertex operator algebra $A$ occurring below, $\Rep(A)$ denotes its category of ordinary modules with the ribbon structure. The following theorem summarizes the comparison.

\begin{theorem}\label{thm:main}
Let $U$ and $V$ be the vertex operator algebras defined in \eqref{eq:Hdef} and \eqref{eq:Jdef}.
\begin{enumerate}
\item The lattice isomorphism in Section~\ref{sub:Riesen} induces an isomorphism $\Psi:U_{\rm Riesen}\longrightarrow U$.
\item There are automorphisms $\tau$ of $U$ and $\sigma$ of $V$ of order two, and an explicit lattice vertex operator algebra isomorphism $\Phi$ on the common lattice subalgebra, such that
\[
 \Phi(U^\tau)=V^\sigma=:C.
\]
The algebra $C$ contains $K(\slc,3)\otimes K(\slc,13)$, where $K(\slc,k)$ is the commutant of the Cartan Heisenberg algebra in $L_{\slc}(k,0)$.
\item The vertex operator algebras $U,V,C$ are simple, unitary, rational and $C_2$-cofinite, and are of CFT type. The algebras $U$ and $V$ are order-two simple current extensions of $C$.
\item There are ribbon equivalences
\[
 \Rep(C)\simeq\Tc\boxtimes\D,
 \qquad \Rep(U)\simeq\D\simeq\Rep(V),
\]
{where $\Tc$ is the full fusion subcategory of $\Rep(C)$ generated by the two odd summands and $\D$ is its M\"uger centralizer.}
\item The associated conformal nets $\mathcal A_U$, $\mathcal A_V$ and $\mathcal A_C$ are completely rational, and
\[
 \mathcal A_U^{\langle\tau\rangle}\cong\mathcal A_C
       \cong\mathcal A_V^{\langle\sigma\rangle}.
\]
\end{enumerate}
\end{theorem}

\begingroup
In part (4), the category $\Tc$ has exactly four simple objects up to isomorphism, namely $\one,\mathsf a,\mathsf b,\mathsf a\mathsf b$. Here $\one=C$ is the tensor unit, while $\mathsf a$ and $\mathsf b$ are the $-1$ eigenspaces of $\tau$ on $U$ and $\sigma$ on $V$, respectively, regarded as $C$-modules through the fixed point isomorphism. Juxtaposition denotes fusion over $C$. All four objects are invertible, so $\Tc$ is pointed. Their fusion group is $(\Z/2)^2$, and their twists are
\[
 \theta_{\one}=\theta_{\mathsf a}=\theta_{\mathsf b}=1,
 \qquad \theta_{\mathsf a\mathsf b}=-1.
\]
For a simple $C$-module $X$ of lowest conformal weight $h_X$, the twist is $\theta_X=e^{2\pi i h_X}$. Thus $\mathsf a$ and $\mathsf b$ have integral conformal weights, while their fusion product has half-integral conformal weight. These data give a nondegenerate subcategory, and $\Tc$ is ribbon equivalent to $Z(\mathrm{Vec}_{\Z/2})$, as shown in Section~\ref{sec:categories}. Here $\mathrm{Vec}_{\Z/2}$ is the category of finite-dimensional $\Z/2$-graded complex vector spaces with its usual associator, $Z$ denotes the Drinfeld center, and $\boxtimes$ between categories denotes the Deligne tensor product. The M\"uger centralizer consists of the objects whose double braiding with every object of $\Tc$ is the identity.
\par\endgroup
An order-two simple current extension is an extension by an irreducible module whose fusion square is the vacuum module. The nets in part (5) and their fixed point subnets are described in Section~\ref{sec:nets}.

Taking a fixed point subalgebra and then a different extension has appeared in other important constructions. Let $\Lambda$ be the Leech lattice, and let $\vartheta$ be the standard lift of its $-1$ isometry to $V_\Lambda$. The fixed point subalgebra is $V_\Lambda^+$. Adjoining its odd untwisted module gives back $V_\Lambda$, whereas adjoining the appropriate integral-weight eigenspace of the $\vartheta$-twisted module gives the moonshine vertex operator algebra $V^\natural$ \cite{FLM,HuangMoonshine}. The two resulting vertex operator algebras are different extensions of the same fixed point subalgebra.

Similarly, the two holomorphic vertex operator algebras of central charge sixteen, $V_{E_8\oplus E_8}$ and $V_{D_{16}^+}$, are distinct order-two simple current extensions of a common fixed point subalgebra. They have the same character but different weight-one Lie algebras \cite{DM}.

The rationality argument also leads to a general result. Huang proved in \cite[Section 4, Proposition 3]{Huang} that his positive conformal-vector identity, together with the embedded regular parafermion subalgebras, implies strong rationality and positivity of the nonvacuum lowest conformal weights. There is a related structure in parafermion vertex operator algebras: their conformal vectors are positive linear combinations of the Virasoro vectors associated to the rank-one parafermion subalgebras \cite[Sections 2--3]{DW}. These subalgebras need not commute. Motivated by Huang's result and this parafermion decomposition, we give a general rationality criterion in Section~\ref{sec:criterion}. Under the stated regularity, conformal-weight and self-adjointness hypotheses, an identity
\[
 \omega_A=\sum_i t_i e_i,\qquad t_i>0,
\]
implies strong rationality of a simple unitary vertex operator algebra $A$; here $e_i$ are the conformal vectors of its specified subalgebras. We include the proof of vacuum projectivity, the $C_2$-cofiniteness argument and the passage from the finite braided category to regularity. This criterion is independent of the Haagerup realization.

For the application, the standard parafermion algebras $K(\slc,3)$ and $K(\slc,13)$ lie in $V^\sigma$. The fixed point isomorphism transports them into $U^\tau$. Averaging their conformal vectors over the compact group generated by the right zero-mode $\slc$ action gives the positive identity needed for the criterion. We compute this identity explicitly in Section~\ref{sec:rationality}. The fixed point comparison itself uses only the lattice generators and a lattice-shift lemma, and is established before rationality is used.

{The two extensions correspond to the algebra objects $\one\oplus\mathsf a$ and $\one\oplus\mathsf b$ in the factor $\Tc$ of $\Rep(C)$. The category of local modules over either algebra in $\Tc$ is equivalent to $\mathrm{Vec}$, the category of finite-dimensional complex vector spaces. The extension theory therefore identifies both $\Rep(U)$ and $\Rep(V)$ with $\D$, as ribbon categories.} An identification of either module category with the Haagerup center consequently gives the same identification for the other. The character, the full weight-one spaces and the identification with the Haagerup center are results from the original constructions; our proof of the common fixed point theorem does not assume them. The unitary lattice embeddings and regularity also give completely rational conformal nets, using the finite-index results for conformal subnets.

The paper is organized as follows. Section~\ref{sec:notation} fixes notation and records the conventions used throughout the paper. Section~\ref{sec:criterion} proves the general rationality criterion. Section~\ref{sec:lattice} gives the lattice construction and the three sets of generators. Sections~\ref{sec:neighbor} and \ref{sec:lift} describe the neighboring lattices and the lift identifying the required generators. Section~\ref{sec:fixed} proves the common fixed point theorem and identifies the parafermion subalgebras. Section~\ref{sec:rationality} proves the positive conformal-vector identity and applies the criterion. Sections~\ref{sec:categories} and \ref{sec:nets} treat the module categories and the conformal nets, respectively. Section~\ref{sec:remaining} discusses further questions. The appendix supplies the complete finite input tables, intermediate expressions, Gram matrices and coefficient identities used in the calculations.

\section{Notation and conventions}\label{sec:notation}
We fix the conventions for vertex operator algebras, modules and their tensor categories. The principal notation for the lattice realizations and their comparison is collected in Table~\ref{tab:notation}.

All vector spaces, algebras and categories are over $\C$, and $\mathbb N=\{0,1,2,\ldots\}$. The symbols $\operatorname{Hom}$ and $\operatorname{End}$ denote spaces of morphisms and endomorphisms in the indicated category. For a vertex operator algebra $A$ with vacuum $\one$ and conformal vector $\omega_A$, we write
\[
 Y_A(a,z)=\sum_{n\in\Z}a_nz^{-n-1},\qquad
 Y_A(\omega_A,z)=\sum_{n\in\Z}L_A(n)z^{-n-2}.
\]
The same mode notation is used on modules; the module is indicated when necessary. We omit the subscript $A$ when there is no ambiguity. The space $A_n$ is the $L_A(0)$-eigenspace of eigenvalue $n$, and $\wt(a)=n$ for $a\in A_n$. CFT type means $A=\bigoplus_{n\ge0}A_n$ and $A_0=\C\one$. The translation operator is $\partial=L_A(-1)$. We use $\nn{ab}=a_{-1}b$; repeated normal products are right associated. Angle brackets with subscript ${\rm VA}$ denote the vertex subalgebra generated using all modes.

A weak module satisfies the vertex operator algebra module identities without an imposed grading. An admissible module has a compatible nonnegative integer grading. An ordinary module has diagonalizable $L(0)$, finite-dimensional weight spaces and weights bounded below in each integer coset. Rationality means complete reducibility of admissible modules; regularity means that every weak module is a direct sum of irreducible ordinary modules. For CFT-type vertex operator algebras, regularity is equivalent to rationality and $C_2$-cofiniteness \cite{ABD}, where
\[
 C_2(A)=\operatorname{span}_{\C}\{a_{-2}b:a,b\in A\}.
\]
Strong rationality means simplicity, CFT type, self-contragredience, rationality and $C_2$-cofiniteness. A rational vertex operator algebra is called holomorphic if its only irreducible module, up to isomorphism, is itself. A unitary vertex operator algebra has an antilinear involution and a positive definite invariant Hermitian form, normalized to give the vacuum norm one \cite{DLUnitary}.

A generalized $A$-module is a weak $A$-module $N$ which is the direct sum of its generalized $L(0)$-eigenspaces:
\[
 N=\bigoplus_{h\in\C}N_{[h]},\qquad
 N_{[h]}=\{v\in N:(L(0)-h)^jv=0\text{ for some }j\ge1\}.
\]
Thus $L(0)$ need not act semisimply. Such a module is called grading-restricted if every $N_{[h]}$ is finite-dimensional and, for each $h\in\C$, the spaces $N_{[h-n]}$ vanish for sufficiently large integers $n$. We write $\mathcal G(A)$ for the category of grading-restricted generalized $A$-modules. The restricted dual $N'=\bigoplus_hN_{[h]}^*$, equipped with the contragredient action, is the contragredient module. Self-contragredience means $A\cong A'$ as $A$-modules. A projective object is one for which the functor $\operatorname{Hom}(P,-)$ is exact.

A vertex operator subalgebra $A$ of $D$ need not have the same conformal vector as $D$. The inclusion $A\subset D$ is called conformal if their conformal vectors agree. We use
\[
 C_D(A)=\{d\in D:a_nd=0\text{ for all }a\in A, n\ge0\}
\]
for the commutant of $A$ in a vertex operator algebra $D$. In the rational setting, $\Rep(A)$ denotes the ribbon category of ordinary modules. Its tensor product is written $\boxtimes$, and a simple current is an irreducible module invertible under this product. For categories, $\boxtimes$ denotes the Deligne tensor product; $\simeq_{\rm rib}$ denotes an equivalence preserving the ribbon structure. The category of finite-dimensional complex vector spaces is denoted by $\mathrm{Vec}$.

The main notation used in this paper is collected in the following table. The indicated sections give the explicit definitions; letters used only inside a lemma or a proof are defined there.
\begingroup\small
\renewcommand{\arraystretch}{1.12}
\begin{longtable}{@{}p{.23\textwidth}p{.71\textwidth}@{}}
\caption{Principal notation}\label{tab:notation}\\
\toprule Symbol & Meaning\\\midrule\endfirsthead
\toprule Symbol & Meaning (continued)\\\midrule\endhead
\bottomrule\endfoot
$B,E,\omega$ & $B=V_{E_6\oplus A_2}\subset E=V_{E_8}$; $\omega$ is their common conformal vector, of central charge $8$.\\
$U_{\rm Riesen},U,V$ & Riesen's, Huang's and Jiao's vertex operator algebras, respectively (Section~\ref{sec:lattice}).\\
$A_{1,k},K(\slc,k)$ & $A_{1,k}=L_{\slc}(k,0)$; $K(\slc,k)$ is the commutant of its Cartan Heisenberg subalgebra.\\
$L,G,\eps$ & The lattice $E_6\oplus A_2$, its Gram matrix and its chosen cocycle (Section~\ref{sec:lattice}).\\
{$\alpha_i^H,\alpha_i^J$} & {Ordered root bases of Huang's original copy $L$ and Jiao's copy $L_J$, respectively; both have Gram matrix $G$.}\\
$r,\xi$ & The constants $r=\sqrt{13}>0$ and $\xi=i\sqrt3$; formal variables are denoted by $z$ or $t$.\\
$\ell,h$ & Huang's level-$1$ and level-$39$ Cartan lattice vectors; $H_L=\ell(-1)\one$, $H_R=h(-1)\one$.\\
$u,u_{2,6},u_{1,7}$ & The additional weight-two generators of $U$ and $V$ (Section~\ref{sec:lattice}).\\
$\omega_L,\omega_R$ & Sugawara conformal vectors of the level-$1$ and level-$39$ affine subalgebras of $U$.\\
$\omega_c,e,f$ & $\omega_c=\omega-\omega_L-\omega_R=e+f$; $e,f$ are the recovered commuting Virasoro vectors.\\
$p,b;\lambda,s,t$ & Jiao's Cartan lattice vectors and $\lambda=(p+b)/2$, $s=(p+b)/4$, $t=(13p-3b)/4$. Here $t$ is a lattice vector.\\
$L_J,L_0,L_H$ & Jiao's coordinate lattice, $L_0=\{x\in L_J:(\lambda,x)\in2\Z\}$ and its even neighbor $L_H=L_0+\Z s$.\\
$\mathsf B$ & The matrix whose columns form the ordered basis of $L_0$ used for the cocycle lifts (Section~\ref{sec:lift}).\\
$\Psi,R_*,\Phi$ & The Riesen--Huang lift, the corrected rational lattice isometry, and its lift for the common fixed point comparison.\\
$\tau,\sigma,C$ & The two involutions and $C=\Phi(U^\tau)=V^\sigma$ (Sections~\ref{sec:neighbor}--\ref{sec:fixed}).\\
$\widetilde U,W$ & The transported copy of $U$ in $V_{L_H}$ and its explicitly generated even subalgebra; $W=C$ by Theorem~\ref{thm:fixed}.\\
$\theta_H,\theta_J$ & The antilinear involutions defining the two unitary structures (Section~\ref{sec:lattice}).\\
$\eta_k^J,\eta_k$ & Parafermion conformal vectors in $V^\sigma$ and their images $\eta_k=\Phi^{-1}(\eta_k^J)\in U^\tau$, for $k=3,13$ (Section~\ref{sec:rationality}).\\
$\mathcal E_R$ & Normalized Haar averaging over the compact group generated by the right zero-mode $\slc$ action (Section~\ref{sec:rationality}).\\
{$M_s,M_p,\Tc,\D$} & {The two odd $C$-modules, the pointed subcategory they generate, with four simple objects, and its M\"uger centralizer (Section~\ref{sec:categories}).}\\
$\mathcal A_U,\mathcal A_V,\mathcal A_C$ & The associated conformal nets (Section~\ref{sec:nets}).\\
\end{longtable}
\endgroup

\section{A rationality criterion}\label{sec:criterion}
In this section, we prove a general strong rationality criterion for a simple unitary vertex operator algebra whose conformal vector is a positive linear combination of conformal vectors of regular subalgebras. We establish projectivity of the vacuum module, prove the categorical implication from projectivity to regularity, and obtain $C_2$-cofiniteness from the positive Hermitian form. The criterion is independent of the Haagerup constructions.

\subsection{The criterion}
The proof of the criterion expands the argument of \cite[Section 4, Proposition 3]{Huang}. The standard tensor-category existence results and the rigidity theorem used below are stated with their precise hypotheses.

\begin{theorem}\label{thm:positivecriterion}
Let $A$ be a simple unitary vertex operator algebra of CFT type with conformal vector $\omega_A$. Suppose that $A$ contains finitely many simple CFT-type vertex operator subalgebras $A^1,\ldots,A^m$, with conformal vectors $e_1,\ldots,e_m$. The embeddings preserve the vacuum, vertex operators and translation. Assume:
\begin{enumerate}
\item $e_i\in A_2$, and each $A^i$ is regular.
\item Every irreducible $A^i$-module has real, nonnegative lowest conformal weight, and the vacuum module $A^i$ is the only irreducible module of lowest weight zero.
\item There are positive real numbers $t_1,\ldots,t_m$ such that
\begin{equation}\label{eq:positivefinite}
 \omega_A=\sum_{i=1}^m t_i e_i.
\end{equation}
\item The operators $L_i(0):=(e_i)_1$ are self-adjoint for the same positive invariant Hermitian form on $A$.
\end{enumerate}
Then $A$ is strongly rational and regular. Every irreducible $A$-module other than $A$ has strictly positive lowest conformal weight.
\end{theorem}

We first establish projectivity of the vacuum module in Lemma~\ref{lem:vacuumprojective} and prove in Proposition~\ref{prop:projectivecriterion} that, under $C_2$-cofiniteness, this implies strong rationality and regularity. The proof of Theorem~\ref{thm:positivecriterion} is completed in Subsection~\ref{subsec:criterionproof}, where $C_2$-cofiniteness is obtained from the positive Hermitian form.

The subalgebras $A^i$ need not commute, and the $t_i$ need not sum to one. No Hermitian form is assumed on an arbitrary $A$-module. Condition (4) holds for unitary embeddings with inherited unitary structures and $L_A(1)e_i=0$, by invariance of the Hermitian form. On any weak module we use the notation
\[
 Y(e_i,z)=\sum_{n\in\Z}L_i(n)z^{-n-2},\qquad
 L_i(n)=(e_i)_{n+1}.
\]
Thus \eqref{eq:positivefinite} implies $L_A(n)=\sum_i t_iL_i(n)$ for every integer $n$.

\subsection{The vacuum module}
\begin{lemma}\label{lem:vacuumprojective}
Let $A$ be a CFT-type vertex operator algebra satisfying the subalgebra hypotheses and conditions (1)--(3) of Theorem~\ref{thm:positivecriterion}. No unitarity assumption is required. For every $N\in\mathcal G(A)$, all generalized conformal weights are real and nonnegative. Moreover, $A$ is projective in $\mathcal G(A)$. If $A$ is simple, then $A\cong A'$, and an irreducible object of $\mathcal G(A)$ with lowest weight zero is isomorphic to $A$.
\end{lemma}
\begin{proof}
Let $N_{[h]}\ne0$. Since $e_i\in A_2$, the commutator formula gives $[L_A(0),L_i(0)]=0$. Hence $L_i(0)$ preserves the finite-dimensional space $N_{[h]}$. The restriction of $N$ to $A^i$ is a weak module and is a direct sum of irreducible ordinary modules by regularity. Condition (2) shows that $L_i(0)$ is diagonalizable with nonnegative real eigenvalues. The same is true of its restriction to $N_{[h]}$. Taking the ordinary linear trace on this space gives
\begin{equation}\label{eq:tracepositivity}
 h\dim N_{[h]}=\sum_i t_i\operatorname{Tr}_{N_{[h]}}L_i(0)\ge0.
\end{equation}
The nilpotent part of $L_A(0)$ has trace zero. It follows that $h$ is real and nonnegative.

For $h=0$, positivity of the $t_i$ forces every eigenvalue of every $L_i(0)|_{N_{[0]}}$ to vanish. Thus $L_i(0)N_{[0]}=0$. By condition (2) and CFT type of $A^i$, the kernel of $L_i(0)$ in $N$ consists of the vacuum lines in copies of $A^i$. Therefore $L_i(-1)N_{[0]}=0$, and \eqref{eq:positivefinite} gives $L_A(-1)N_{[0]}=0$.

It follows that each $v\in N_{[0]}$ is a vacuum-like vector, that is, $a_nv=0$ for all $a\in A$ and $n\ge0$. By \cite[Propositions 3.3--3.4]{LiForms}, $v$ therefore determines an $A$-module homomorphism $s_v:A\to N$ satisfying $s_v(\one)=v$.

To prove that $A$ is projective in $\mathcal G(A)$, it suffices to show that every epimorphism onto $A$ splits. Let $\pi:N\twoheadrightarrow A$ be any such epimorphism. Since module maps preserve generalized weights, there is $v\in N_{[0]}$ with $\pi(v)=\one$. The composition $\pi s_v$ fixes $\one$ and is therefore the identity on $A$. Thus $\pi$ splits. This proves projectivity. If $N$ is irreducible with lowest weight zero, any nonzero $v\in N_{[0]}$ gives a nonzero map $s_v$, which is an isomorphism when $A$ is simple.

Finally, take $N=A'$. Its weight-zero space is one-dimensional. A nonzero vector in this space therefore determines a nonzero module homomorphism $A\to A'$ by the same result of Li. For simple $A$ this map is injective. Since $A$ and $A'$ have the same finite dimension in each weight, the map is also surjective. Thus $A\cong A'$ as $A$-modules.
\end{proof}

\subsection{The finite braided category}
We need a growth estimate before rigidity is known. A finite $\C$-linear abelian category means one equivalent to the category of finite-dimensional modules over a finite-dimensional associative $\C$-algebra. In particular it has finite-dimensional morphism spaces, finitely many simple objects and projective covers.

\begin{lemma}\label{lem:moderategrowth}
Let $\mathcal C$ be a finite $\C$-linear abelian monoidal category with tensor product $\boxtimes$, unit $\mathbf1_{\mathcal C}$ and right-exact tensor product in each variable. For every object $X$, there are constants $a,b\ge1$ such that
\[
 \dim\operatorname{End}_{\mathcal C}(X^{\boxtimes n})\le b a^{2n}
 \quad(n\ge0).
\]
In particular, this dimension is less than $n!$ for sufficiently large $n$.
\end{lemma}
\begin{proof}
Let $P$ be the direct sum of projective covers of representatives of the finitely many simple objects. It is a projective generator: every object is a quotient of a finite direct sum of copies of $P$. Put $d(Y)=\dim\operatorname{Hom}_{\mathcal C}(P,Y)$. The functor $\operatorname{Hom}_{\mathcal C}(P,-)$ is exact and faithful. A basis of $\operatorname{Hom}_{\mathcal C}(P,Y)$ gives an epimorphism $P^{\oplus d(Y)}\twoheadrightarrow Y$. Applying the right-exact functor $X\boxtimes-$ and then the exact functor $\operatorname{Hom}_{\mathcal C}(P,-)$ gives
\[
 d(X\boxtimes Y)\le d(X\boxtimes P)d(Y).
\]
Set $a=\max\{1,d(X\boxtimes P)\}$. Starting with $X^{\boxtimes0}=\mathbf1_{\mathcal C}$, induction gives $d(X^{\boxtimes n})\le a^n d(\mathbf1_{\mathcal C})$. Faithfulness gives an injection
\[
 \operatorname{End}_{\mathcal C}(Y)\hookrightarrow
 \operatorname{End}_{\C}\bigl(\operatorname{Hom}_{\mathcal C}(P,Y)\bigr).
\]
The assertion follows with $b=\max\{1,d(\mathbf1_{\mathcal C})^2\}$, since factorial growth eventually exceeds every fixed exponential.
\end{proof}

\begin{proposition}\label{prop:projectivecriterion}
Let $A$ be simple, self-contragredient, $C_2$-cofinite and of CFT type. If its vacuum module is projective in $\mathcal G(A)$, then $A$ is strongly rational and regular.
\end{proposition}
\begin{proof}
The standard vertex tensor-category results make $\mathcal C=\mathcal G(A)$ a finite $\C$-linear abelian braided monoidal category, with unit $A$ and tensor product right exact in each variable. Contragredient duality and $A\cong A'$ give natural isomorphisms
\begin{equation}\label{eq:GVduality}
 \operatorname{Hom}_{\mathcal C}(Y\boxtimes X,A)
 \cong\operatorname{Hom}_{\mathcal C}(Y,X').
\end{equation}
These results do not assume rigidity; see \cite[Sections 2.3 and 3]{McRae}. In particular, the identity of $X'$ corresponds to a nonzero evaluation morphism $\mathrm{ev}_X:X'\boxtimes X\to A$ for every nonzero $X$. Since $A$ is simple, this morphism is surjective. Since $A$ is projective, it splits. Thus $A$ is a direct summand of $X'\boxtimes X$.

We use \cite[Theorem 1.1]{EP}: in a $\C$-linear braided category with finite-dimensional morphism spaces, finite direct sums, split idempotents and $\operatorname{End}(\mathbf1)=\C$, an indecomposable object $X$ is rigid if the unit is a direct summand of $Y\boxtimes X$ for some $Y$ and $\dim\operatorname{End}(X^{\boxtimes n})<n!$ for some $n$. These two conditions are called non-negligibility and moderate growth, respectively. All the categorical hypotheses hold for $\mathcal C$; in particular $\operatorname{End}_A(A)=\C$. The splitting above verifies non-negligibility, and Lemma~\ref{lem:moderategrowth} verifies moderate growth. Every indecomposable object is therefore rigid. Each object is a finite direct sum of indecomposable objects, so $\mathcal C$ is rigid.

Tensoring in a rigid category has both left and right adjoints and is consequently exact. For a dual $X^\vee$ chosen with the corresponding adjunction, there is a natural isomorphism
\[
 \operatorname{Hom}_{\mathcal C}(X,-)
 \cong\operatorname{Hom}_{\mathcal C}(A,-\boxtimes X^\vee).
\]
The right-hand functor is exact because $A$ is projective. Hence every object is projective and every short exact sequence splits. Thus $\mathcal C$ is semisimple.

It remains to pass to arbitrary weak modules. For an $\mathbb N$-graded $C_2$-cofinite vertex operator algebra, the finitely generated weak modules are exactly the grading-restricted generalized modules \cite[Proposition 2.11]{McRae}. Let $N$ be any weak $A$-module and $v\in N$. Its cyclic submodule $\langle v\rangle_A$, generated by all modes of $A$, belongs to $\mathcal C$ and is therefore a finite direct sum of irreducible modules. Such irreducible generalized modules are ordinary: the locally nilpotent part of $L(0)$ commutes with all modes, and its kernel is a nonzero submodule, so simplicity makes that part zero. Consequently $N$ is a sum of irreducible ordinary submodules. A module which is a sum of simple submodules is semisimple, and hence a direct sum of them. This proves regularity and therefore rationality. All the other conditions for strong rationality were assumed.
\end{proof}

\subsection{Proof of the criterion}\label{subsec:criterionproof}
\begin{proof}[Proof of Theorem~\ref{thm:positivecriterion}]
We prove $C_2$-cofiniteness using the positive Hermitian form, denoted by $(\cdot,\cdot)_A$, with $\|v\|^2=(v,v)_A$. For an $A^i$-module $M$, define
\[
 C_2^{A^i}(M)=\operatorname{span}_{\C}\{a_{-2}m:a\in A^i,\ m\in M\}.
\]
Regularity and CFT type imply rationality and $C_2$-cofiniteness of $A^i$ \cite[Theorem 4.5]{ABD}. It has finitely many irreducible ordinary modules, and $M/C_2^{A^i}(M)$ is finite-dimensional for each such module \cite[Proposition 5.2]{ABD}. This quotient is graded by $L_i(0)$. Thus there is a real number $D_i\ge0$, uniform over all these irreducible modules, such that
\[
 M(\lambda)\subset C_2^{A^i}(M)\quad(\lambda>D_i),
 \qquad M(\lambda)=\ker_M(L_i(0)-\lambda).
\]
Here $\lambda$ is the actual eigenvalue, including the lowest conformal weight.

Restricting $A$ to $A^i$ gives a direct sum of these irreducible modules. Every vector belongs to finitely many summands. The same bound therefore gives
\begin{equation}\label{eq:highpartialweights}
 \bigoplus_{\lambda>D_i}\ker_A(L_i(0)-\lambda)
 \subset C_2^{A^i}(A)\subset C_2(A).
\end{equation}
The bound does not depend on the multiplicities in the restriction.

For $n\ge0$, take the orthogonal complement inside the finite-dimensional space $A_n$:
\[
 W_n=(C_2(A)\cap A_n)^\perp\subset A_n.
\]
Each $L_i(0)$ preserves $A_n$ and is self-adjoint. Its eigenspaces with eigenvalues greater than $D_i$ lie in $C_2(A)\cap A_n$ by \eqref{eq:highpartialweights}. Hence a vector $v\in W_n$ has zero projection onto those eigenspaces. Applying the spectral theorem separately to each self-adjoint operator $L_i(0)$ on $A_n$, we obtain $(L_i(0)v,v)_A\le D_i\|v\|^2$. Summing with the positive coefficients yields
\[
 n\|v\|^2=\sum_i t_i(L_i(0)v,v)_A
 \le D\|v\|^2,\qquad D=\sum_i t_iD_i.
\]
Thus $W_n=0$ for $n>D$. Since $C_2(A)$ is a homogeneous subspace of $A$,
\[
 \dim A/C_2(A)\le\sum_{0\le n\le\lfloor D\rfloor}\dim A_n<\infty.
\]

Lemma~\ref{lem:vacuumprojective} gives self-contragredience and projectivity of the vacuum module. Proposition~\ref{prop:projectivecriterion} now gives strong rationality and regularity. Finally \eqref{eq:tracepositivity} shows that every lowest conformal weight is nonnegative, and Lemma~\ref{lem:vacuumprojective} identifies the only irreducible module with lowest weight zero as $A$.
\end{proof}

\begin{corollary}\label{cor:nonunitarycriterion}
Let $A$ be a simple CFT-type vertex operator algebra satisfying the subalgebra hypotheses and conditions (1)--(3) of Theorem~\ref{thm:positivecriterion}. If $A$ is $C_2$-cofinite, then it is strongly rational and regular, and its nonvacuum irreducible modules have strictly positive lowest conformal weights. No unitarity assumption is needed.
\end{corollary}
\begin{proof}
Lemma~\ref{lem:vacuumprojective} supplies self-contragredience, vacuum projectivity and the assertion about conformal weights. Proposition~\ref{prop:projectivecriterion} applies.
\end{proof}

\begin{remark}\label{rem:paramotivation}
The parafermion decomposition mentioned in the Introduction has a precise form. Let $\mathfrak g$ be simple with positive roots $\Delta_+$ and dual Coxeter number $h^\vee$, and normalize long roots to have squared length $2$. For a positive integer $k$, put $k_\alpha=2k/(\alpha,\alpha)$. The root-$\alpha$ subalgebra $K(\slc,k_\alpha)$ has Virasoro vector $\omega_\alpha$, and the conformal vector of $K(\mathfrak g,k)$ satisfies
\[
 \omega_{K(\mathfrak g,k)}
 =\sum_{\alpha\in\Delta_+}
   \frac{k(k_\alpha+2)}{k_\alpha(k+h^\vee)}\,\omega_\alpha
 \qquad\text{\cite[Section 2]{DW}}.
\]
Every coefficient is positive. Here $K(\mathfrak g,k)$ is the commutant of the Cartan Heisenberg subalgebra in $L_{\mathfrak g}(k,0)$. The positivity and self-adjointness used in \cite[Section 3]{DW} provide one motivation for the criterion. Regularity is required of the parafermion subalgebras, not of the Virasoro algebras generated by their conformal vectors alone.
\end{remark}

\section{Lattice vertex operator algebras and generators of the Haagerup realization VOAs}\label{sec:lattice}
We describe the inclusion of the lattice $E_6\oplus A_2$ in $E_8$ and fix the ordered root basis, cocycle and unitary structures. The generators of the three Haagerup realizations are then expressed in these conventions. We recover Riesen's two commuting Virasoro vectors from Huang's generators and prove that the Riesen and Huang constructions are isomorphic.

\subsection{The lattices \texorpdfstring{$E_6\oplus A_2$ and $E_8$}{E6 + A2 and E8}}\label{subsec:latticeconventions}
Let $L=\bigoplus_{i=1}^8\Z\alpha_i$ be the root lattice $E_6\oplus A_2$. The entries of its Gram matrix $G=((\alpha_i,\alpha_j))_{1\le i,j\le8}$ are
\begin{equation}\label{eq:grammatrix}
 G_{ij}=\begin{cases}
 2,&i=j,\\
 -1,&(i,j)\in\mathcal P\text{ or }(j,i)\in\mathcal P,\\
 0,&\text{otherwise},
 \end{cases}
 \qquad
 \mathcal P=\{(1,2),(2,3),(3,4),(4,5),(3,6),(7,8)\}.
\end{equation}
Thus $\alpha_1,\ldots,\alpha_6$ are simple roots of the $E_6$ factor and $\alpha_7,\alpha_8$ are simple roots of the $A_2$ factor. We realize these eight roots in an $E_8$ lattice by adjoining the glue vector below; they form a basis of $L$, not a simple-root basis of $E_8$. All lattice coordinates are taken in the indicated ordered basis: $x=(x_1,\ldots,x_8)^T$ means $x=\sum_{i=1}^8x_i\alpha_i$. Unless a different range is specified, lattice-coordinate indices run from $1$ to $8$. The bilinear form is $(x,y)=x^TGy$, and $x^2=(x,x)$. We choose the bimultiplicative cocycle $\eps:L\times L\to\{\pm1\}$ determined by
\[
 \eps(\alpha_i,\alpha_j)=\begin{cases}
 (-1)^{(\alpha_i,\alpha_j)},&i>j,\\
 -1,&i=j,\\
 1,&i<j.
 \end{cases}
\]
Thus, for $x=\sum_i x_i\alpha_i$ and $y=\sum_i y_i\alpha_i$, we have
\begin{equation}\label{eq:cocycle}
 \eps(x,y)=(-1)^{\sum_i x_i y_i+\sum_{i>j}x_iG_{ij}y_j}.
\end{equation}
The cocycle $\eps$ satisfies
\[
 \eps(x,y)\eps(y,x)^{-1}=(-1)^{(x,y)},\qquad
 \eps(x,x)=(-1)^{x^2/2}.
\]
In particular, for a root $\alpha$, the standard $\slc$-triple is
\[
 \bigl(\alpha(-1)\one,e^\alpha,-e^{-\alpha}\bigr).
\]
The minus sign in the last entry is part of the convention.

Set $B=V_L=M(1)\otimes\C_\eps[L]$. Here $\mathfrak h=\C\otimes_\Z L$, $M(1)$ is its Heisenberg vacuum module, and $\C_\eps[L]$ is the twisted group algebra with basis $e^x$ and multiplication $e^xe^y=\eps(x,y)e^{x+y}$. We also write $e^x$ for $\one\otimes e^x$. For $h\in\mathfrak h$, $h(n)$ is the $n$th Heisenberg mode, so $h(n)$ is the mode $(h(-1)\one)_n$. We use
\[
 Y(a,z)=\sum_{n\in\Z}a_{n}z^{-n-1},\qquad
 Y(e^x,z)e^y=\eps(x,y)z^{(x,y)}
 \exp\!\left(\sum_{n\ge1}\frac{x(-n)}n z^n\right)e^{x+y}.
\]
The lattice conformal vector is
\begin{equation}\label{eq:omega}
 \omega=\frac12\sum_{i,j}(G^{-1})_{ij}\alpha_i(-1)\alpha_j(-1)\one,
 \qquad c=8.
\end{equation}
Thus $\wt(e^x)=x^2/2$. The lattice vertex operator algebra is generated by the eight Heisenberg elements and $e^{\pm\alpha_i}$. Throughout this paper, generation means generation as a vertex algebra using all modes.

The inclusion in $E_8$ is equally explicit. Put
\begin{equation}\label{eq:glue}
 \mu=(4/3,5/3,2,4/3,2/3,1,2/3,1/3)^T,\qquad \widetilde L=L+\Z\mu.
\end{equation}
One checks that $G\mu=(1,0,0,0,0,0,1,0)^T$, $\mu^2=2$, $3\mu\in L$, and $\det G=9$. Hence $\widetilde L$ is even unimodular of rank eight and is isometric to $E_8$. A compatible cocycle extension gives
\begin{equation}
 V_{E_8}=V_L\oplus V_{L+\mu}\oplus V_{L-\mu}.
\end{equation}
The order-three grading automorphism, acting on the three summands by $1,\zeta_3,\zeta_3^2$, has fixed point subalgebra $B$, where $\zeta_3=e^{2\pi i/3}$. This order-three construction is distinct from the order-two comparison used later.

\subsection{Huang's level-\texorpdfstring{$(1,39)$}{(1,39)} generators}
Using the notation of Subsection~\ref{subsec:latticeconventions}, we give the generators of Huang's vertex operator algebra. Define
\[
 \ell=(-1,-2,-3,-2,-1,-2,0,0)^T,\qquad
 h=(5,8,9,8,5,0,2,2)^T,\qquad I=\{1,2,3,4,5,7,8\}.
\]
Let $c_i=h_i$ be the $i$th coordinate of $h$. The six weight-one generators are
\begin{align}
 H_L&=\ell(-1)\one,& E_L&=e^\ell,& F_L&=-e^{-\ell},\label{eq:Hleft}\\
 H_R&=h(-1)\one,& E_R&=\sum_{i\in I}e^{\alpha_i},&
 F_R&=-\sum_{i\in I}c_i e^{-\alpha_i}.\label{eq:Hright}
\end{align}
Indeed $\ell^2=2$, $h^2=78$, and $(\ell,h)=0$. The principal $\slc$-triples in the $A_5$ chain $\alpha_1,\ldots,\alpha_5$ and the external $A_2$ generate affine $\slc$ subalgebras of levels $35$ and $4$, respectively. Their diagonal triple $(H_R,E_R,F_R)$ therefore has level $39$, while the left triple has level $1$. Thus the six elements generate $A_{1,1}\otimes A_{1,39}$, where $A_{1,k}=L_{\slc}(k,0)$.

Throughout set $r=\sqrt{13}>0$ and $\xi =i\sqrt3$. Define
\[
 a=\alpha_7+\alpha_8-\alpha_6,\quad
 d=\alpha_4+\alpha_5-\alpha_6,\quad
 c=\alpha_1+\alpha_2-\alpha_6.
\]
Define the weight-two vector $u$ and the vertex operator algebra $U$ by
\begin{align}
 u&=e^a+\frac{r+3\xi }{20}e^d+\frac{r-3\xi }{20}e^c,\label{eq:Hgenerator}\\
 U&=\langle H_L,E_L,F_L,H_R,E_R,F_R,\omega,u\rangle_{\rm VA}.
 \label{eq:Hdef}
\end{align}
The vector $u$ is a highest weight vector of weight $(1,13)$ for the action of $\slc\oplus\slc$ by zero modes. These formulas reproduce \cite[Section 2, (3)--(8)]{Huang} with its lattice vector $b$ denoted here by $d$. The two Cartan generators can also be obtained from $H_L=(E_L)_{0}F_L$ and $H_R=(E_R)_{0}F_R$.

For a level-$k$ $\slc$-triple $(E,F,H)$ write
\begin{equation}\label{eq:sug}
 \Omega_k(E,F,H)=\frac{\frac12H_{-1}H+E_{-1}F+F_{-1}E}{2(k+2)}.
\end{equation}
Then $\omega_L=\Omega_1(E_L,F_L,H_L)$ and
$\omega_R=\Omega_{39}(E_R,F_R,H_R)$ have central charges $1$ and $117/41$. The remaining conformal vector has central charge $170/41$.

\subsection{Jiao's level-\texorpdfstring{$(3,13)$}{(3,13)} generators in the same root convention}
We use the ordered root basis and the cocycle fixed in Section~\ref{subsec:latticeconventions}. In these conventions, we give the nine generators of Jiao's construction \cite{Jiao} as follows. 

First we put
\begin{align*}
 \rho_1&=(1,1,1,0,0,0,0,0)^T,&
 \rho_2&=(0,1,1,1,0,0,0,0)^T,&
 \rho_3&=(0,0,1,1,1,0,0,0)^T,\\
 q&=(1,2,3,2,1,2,0,0)^T,&
 p&=(1,2,3,2,1,0,0,0)^T,&
 b&=(3,4,3,4,3,2,2,2)^T.
\end{align*}
The six weight-one generators are
\begin{align}
 P_E&=e^{\rho_1}+e^{\rho_2}+e^{\rho_3},\notag\\
 P_F&=-e^{-\rho_1}-e^{-\rho_2}-e^{-\rho_3},\qquad P_H=p(-1)\one,
 \label{eq:P}\\
 B_E&=2e^{\alpha_1}+e^{\alpha_2}+2e^{\alpha_4}+e^{\alpha_5}
       +e^q+2e^{\alpha_7}+e^{\alpha_8},\label{eq:BE}\\
 B_F&=-e^{-\alpha_1}-2e^{-\alpha_2}-e^{-\alpha_4}-2e^{-\alpha_5}
       -e^{-q}-e^{-\alpha_7}-2e^{-\alpha_8},\label{eq:BF}\\
 B_H&=b(-1)\one.\label{eq:BH}
\end{align}
For either $\slc$-triple one has $E_{0}F=H$, $H_{0}E=2E$, $H_{0}F=-2F$, $E_{1}F=k\one$ and $H_{1}H=2k\one$, with $k=3$ and $13$ respectively. The two $\slc$-triples commute. The norms are
\begin{equation}\label{eq:pb}
 p^2=6,\qquad b^2=26,\qquad (p,b)=0.
\end{equation}
Thus they generate $A_{1,3}\otimes A_{1,13}$. The level-three factor is the diagonal sum of the three commuting root
$\slc$-triples associated with the mutually orthogonal roots $\rho_1,\rho_2,\rho_3$, each contributing level $1$. The level-thirteen factor is the diagonal sum of the principal
$\mathfrak{sl}_2$-triples in the three mutually orthogonal $A_2$ root
subsystems generated by $\{\alpha_1,\alpha_2\}$,
$\{\alpha_4,\alpha_5\}$, and $\{\alpha_7,\alpha_8\}$, together with the commuting root triple associated with $q$. These four triples contribute levels $4,4,4,1$, respectively,
giving total level $13$.

Jiao's construction has two additional generators
$u_{1,7},u_{2,6}\in (V_L)_2$.
For $(a,b)=(1,7)$ or $(2,6)$, the vector $u_{a,b}$ is a
highest weight vector of weight $(a,b)$ for the zero-mode action
of the two $\mathfrak{sl}_2$ factors, corresponding respectively
to affine levels $3$ and $13$.
With the cocycle fixed in \eqref{eq:cocycle}, these generators
are given by
\[
 u_{a,b}=\sum_{x\in S_{a,b}}c_{a,b}(x)e^x,
 \qquad (a,b)\in\{(1,7),(2,6)\},
\]
where the supports $S_{a,b}$ and coefficients $c_{a,b}(x)$
are listed in the following table.
Every lattice vector occurring in the table has squared norm
four, so each term $e^x$ has conformal weight two.

\begin{center}\small
\begin{tabular}{cll}\toprule
Generator & Lattice vector $x\in S_{a,b}$ & $c_{a,b}(x)$\\\midrule
$u_{2,6}$ & $(0,1,1,1,1,0,1,1)^T$ & $(r+\xi )/6$\\
& $(1,1,1,1,0,0,1,1)^T$ & $(r-\xi )/3$\\
& $(1,1,1,1,1,0,0,1)^T$ & $(1-\xi )/6$\\
& $(1,1,1,1,1,0,1,0)^T$ & $(1+\xi )/3$\\
& $(1,1,1,2,1,0,0,0)^T$ & $(1-r\xi )/6$\\
& $(1,2,1,1,1,0,0,0)^T$ & $(1+r\xi )/12$\\
& $(2,3,4,3,2,2,0,0)^T$ & $1$\\\midrule
$u_{1,7}$ & $(0,0,-1,0,0,-1,1,1)^T$ & $\xi$\\
& $(1,1,2,2,1,1,1,1)^T$ & $r/2+\xi /6$\\
& $(1,2,2,1,1,1,1,1)^T$ & $r/4-\xi /12$\\
& $(1,2,2,2,1,1,0,1)^T$ & $1/2+\xi /6$\\
& $(1,2,2,2,1,1,1,0)^T$ & $1-\xi /3$\\
& $(1,2,2,2,2,1,0,0)^T$ & $1/2-r\xi /6$\\
& $(2,2,2,2,1,1,0,0)^T$ & $1+r\xi /3$\\\bottomrule
\end{tabular}
\end{center}
Define
\begin{equation}\label{eq:Jdef}
 V=\langle P_E,P_F,P_H,B_E,B_F,B_H,\omega,u_{2,6},u_{1,7}\rangle_{\rm VA}.
\end{equation}
This is the nine-element presentation of \cite[Section 3]{Jiao}. Set $\omega_P=\Omega_3(P_E,P_F,P_H)$ and $\omega_B^{\rm aff}=\Omega_{13}(B_E,B_F,B_H)$. Their central charges are $9/5$ and $13/5$, leaving $18/5$ for the conformal vector of the affine commutant.

\subsection{Unitary structures compatible with the generators}
We now choose positive definite invariant Hermitian forms compatible with the generators in Sections~4.2 and~4.3. These forms will be used to interpret the Casimir projections in Section~\ref{sub:Riesen} as orthogonal projections and to justify the compact-group averages used later. The standard antilinear involution of the lattice vertex operator algebra $V_L$ does not preserve all the displayed $\slc$-subalgebras, since in some triples $(E,F,H)$ the image of $E$ is not a scalar multiple of $F$. We therefore choose a unitary structure for each realization whose antilinear involution preserves its two $\slc$-subalgebras.

Write $(\cdot,\cdot)_0$ and $\theta_0$ for the standard positive definite Hermitian form and antilinear involution of $V_L$; the corresponding unitary structure is given in \cite[Theorem 4.12]{DLUnitary}. A positive lattice character is a group homomorphism $\kappa:L\longrightarrow\mathbb R_{>0}$. Such a character defines an automorphism $g_\kappa$ of $V_L$ by
\[
 g_\kappa(e^x)=\kappa(x)^{-1/2}e^x,
\]
with $g_\kappa$ acting identically on the Heisenberg generators. Set
\[
 (v,w)_\kappa=(g_\kappa^{-1}v,g_\kappa^{-1}w)_0,
 \qquad \theta_\kappa=g_\kappa\theta_0g_\kappa^{-1}.
\]
This gives another unitary structure on the same lattice vertex operator algebra. In particular, for $h_j\in\mathbb R\otimes_{\Z}L$ and integers $n_j>0$, its antilinear involution is
\[
 \theta_\kappa\bigl(h_1(-n_1)\cdots h_m(-n_m)e^x\bigr)
 =(-1)^m\kappa(x)h_1(-n_1)\cdots h_m(-n_m)e^{-x}.
\]

For the generators of Huang and Jiao, respectively, we take
\[
 \kappa_H(x)=5^{x_1}8^{x_2}9^{x_3}8^{x_4}5^{x_5}2^{x_7}2^{x_8},
 \qquad
 \kappa_J(x)=2^{-x_1+x_2-x_4+x_5-x_7+x_8},
\]
where $x=\sum_{i=1}^8x_i\alpha_i$, and denote the corresponding involutions by $\theta_H$ and $\theta_J$. Put $a_H=\kappa_H(\ell)=5^{-2}8^{-4}9^{-3}=1/74649600$. Direct substitution gives
\[
 \begin{aligned}
 \theta_H(E_L)&=-a_HF_L,&\theta_H(F_L)&=-a_H^{-1}E_L,&
 \theta_H(E_R)&=-F_R,\\
 \theta_J(P_E)&=-P_F,&\theta_J(B_E)&=-B_F.
 \end{aligned}
\]
Each of the four Cartan generators is sent to its negative. Consequently, with $\dagger$ denoting the adjoint for the corresponding Hermitian form,
\[
 (E_L)_0^\dagger=a_H(F_L)_0,\qquad
 (E_R)_0^\dagger=(F_R)_0,\qquad
 (P_E)_0^\dagger=(P_F)_0,\qquad
 (B_E)_0^\dagger=(B_F)_0,
\]
and each Cartan zero mode is self-adjoint. For Huang's left triple $(E,F,H)$, we have $E_0^\dagger=a_HF_0$ and $F_0^\dagger=a_H^{-1}E_0$. Taking adjoints reverses the order of a product, so
\[
 (E_0F_0)^\dagger=(a_H^{-1}E_0)(a_HF_0)=E_0F_0,
 \qquad
 (F_0E_0)^\dagger=(a_HF_0)(a_H^{-1}E_0)=F_0E_0.
\]
For the other three triples, $E_0^\dagger=F_0$ and $F_0^\dagger=E_0$, so these identities also hold. Since $H_0^\dagger=H_0$, all four Casimir operators $H_0^2+2E_0F_0+2F_0E_0$ are self-adjoint on each homogeneous subspace of $V_L$ for the corresponding Hermitian forms.

To restrict these unitary structures to $U$ and $V$, one must also check stability under the involutions for the additional weight-two generators. This is established by the lowering identities in and following~\eqref{eq:huangpct}, and the resulting unitary structures on $U$ and $V$ are stated in Corollary~\ref{cor:unitary}. A separate compatibility statement, $\Phi\theta_H=\theta_J\Phi$ for the comparison isomorphism on the common lattice vertex operator algebra, is proved in Proposition~\ref{prop:dictionary}.

\subsection{The generators of Riesen and Huang}\label{sub:Riesen}
We give Riesen's generators in his root coordinates and compare them with Huang's generators under an explicit lattice vertex operator algebra isomorphism. The six affine generators and the additional highest weight vector correspond directly. We then recover the images of Riesen's two Virasoro generators from Huang's generators, which proves equality of the two generated algebras under this identification.

Let $\delta_1,\ldots,\delta_8$ be Riesen's ordered simple-root basis of the root lattice $L_{\rm Rie}\cong A_2\oplus E_6$, renumbered from $1$ to $8$. The first two roots form a basis of the $A_2$ factor and the last six form a basis of the $E_6$ factor. Thus $L_{\rm Rie}=\bigoplus_{i=1}^8\Z\delta_i$ has determinant $9$ and embeds in $E_8$ with index $3$. Write $G_{\rm Rie}=((\delta_i,\delta_j))_{1\le i,j\le8}$ for its Gram matrix. Its diagonal entries are $2$, its off-diagonal entries are $-1$ at the pairs
$(1,2),(3,5),(4,6),(5,6),(6,7),(7,8)$ and their transposes, and all other entries are zero.
Denote by $\eps_{\rm Rie}$ the bimultiplicative cocycle defined by~\eqref{eq:cocycle} in this ordered basis with Gram matrix $G_{\rm Rie}$. This is the cocycle in \cite[Appendix A.1]{RiesenLong}, after increasing each root index by one.

Following \cite[Appendix A.2]{RiesenLong}, put
\[
 (\zeta_1,\ldots,\zeta_7)
 =(\delta_1,\delta_2,\delta_4,\delta_6,\delta_7,\delta_8,
   \delta_3+\delta_4+2\delta_5+2\delta_6+\delta_7),
\]
\[
 (m_1,\ldots,m_7)=(2,2,5,8,9,8,5),\qquad
 h_{\rm Rie}=\sum_{i=1}^7m_i\zeta_i=(2,2,5,10,10,18,14,8)^T.
\]
The coordinates of $h_{\rm Rie}$ are relative to the $\delta_i$. Riesen's six weight-one generators are
\begin{align*}
 H_L^{\rm Rie}&=\delta_3(-1)\one,&
 E_L^{\rm Rie}&=e^{\delta_3},&
 F_L^{\rm Rie}&=-e^{-\delta_3},\\
 H_R^{\rm Rie}&=h_{\rm Rie}(-1)\one,&
 E_R^{\rm Rie}&=\sum_{i=1}^7e^{\zeta_i},&
 F_R^{\rm Rie}&=-\sum_{i=1}^7m_i e^{-\zeta_i}.
\end{align*}
They form two commuting affine $\slc$-triples of levels $1$ and $39$. Let $\omega_L^{\rm Rie}$ and $\omega_R^{\rm Rie}$ be their Sugawara vectors, defined by~\eqref{eq:sug}, and let $\omega_{\rm Rie}$ be the conformal vector of $V_{L_{\rm Rie}}$.

The additional highest weight vector is
\[
 v=\frac{-r+3\xi}{20}e^{\beta_1}
   +\frac{r+3\xi}{20}e^{\beta_2}+e^{\beta_3},
\]
where, in the $\delta_i$ coordinates,
\[
 \beta_1=(0,0,1,2,1,3,2,1)^T,\quad
 \beta_2=(0,0,2,2,3,4,3,2)^T,\quad
 \beta_3=(1,1,1,1,1,2,2,1)^T.
\]
Each $\beta_i$ has squared norm four. This is the vector in \cite[Definition 2.1 and Appendix E]{RiesenLong}, with $r=\sqrt{13}$ and $\xi=i\sqrt3$ as above.

We now give Riesen's two additional Virasoro generators $e_R$ and $f_R$. To express them, we first define four auxiliary weight-two vectors $b_0,b_1,b_2,b_3\in(V_{L_{\rm Rie}})_2$. Their images under the isomorphism $\Psi:V_{L_{\rm Rie}}\longrightarrow V_L$ constructed below are denoted by $B_k=\Psi(b_k)$, $0\le k\le3$. Tables~\ref{tab:Bquadratic} and~\ref{tab:Bexponential} in Appendix~\ref{app:invariants} give the coefficients of the $B_k$ in Huang's coordinates. We use the same coefficient arrays to describe the $b_k$ in Riesen's coordinates, with the cocycle correction specified below. Write $z=\sum_{i=1}^7n_i\zeta_i$ with $n=(n_1,\ldots,n_7)^T\in\Z^7$, and put
\[
 \eta(n)=(-1)^{n_7(n_3+n_4+n_5+n_6)},\qquad
 \mathsf S_{\rm Rie}(z)
 =\eta(n)\left(e^z+\left(\prod_{i=1}^7m_i^{n_i}\right)e^{-z}\right).
\]
For a row $\Gamma$ of Table~\ref{tab:Bexponential}, let $\Gamma_{\rm Rie}$ be the set obtained by replacing every symbol $\nu_i$ in that row by $\zeta_i$. Then
\[
 b_k=\sum_{1\le i\le j\le7}q_{ij}^{(k)}\zeta_i(-1)\zeta_j(-1)\one
       +\sum_{\Gamma}g_\Gamma^{(k)}
          \sum_{z\in\Gamma_{\rm Rie}}\mathsf S_{\rm Rie}(z),
 \qquad 0\le k\le3.
\]
All omitted coefficients are zero. These are the vectors of \cite[Appendix C]{RiesenLong}, with the terms involving $e^z$ and $e^{-z}$ grouped using $\mathsf S_{\rm Rie}(z)$. The factor $\eta(n)$ is the cocycle correction for the seven-root coordinates in that appendix; it is included explicitly here.

Following \cite[Section 2.2, (2.5)--(2.6)]{RiesenLong}, the two Virasoro generators are
\[
 e_R=\frac{2b_0+6b_1}{31}
       +\frac{2-r}{93}b_2+\frac{2-2r}{93}b_3,
 \qquad f_R=\omega_C^{\rm Rie}-e_R,
\]
where
\[
 \omega_C^{\rm Rie}
 =\omega_{\rm Rie}-\omega_L^{\rm Rie}-\omega_R^{\rm Rie}
 =\frac{8b_0+9b_1+2b_2+2b_3}{41}.
\]
Riesen's vertex operator algebra is therefore the explicitly defined nine-generator algebra
\[
 U_{\rm Riesen}
 =\bigl\langle H_L^{\rm Rie},E_L^{\rm Rie},F_L^{\rm Rie},
               H_R^{\rm Rie},E_R^{\rm Rie},F_R^{\rm Rie},
               v,e_R,f_R\bigr\rangle_{\rm VA}.
\]

Define $P:L_{\rm Rie}\longrightarrow L$ by the following images of the basis vectors, expressed in the basis $\alpha_1,\ldots,\alpha_8$ fixed in Section~\ref{subsec:latticeconventions}:
\begin{center}\small
\begin{tabular}{cl}\toprule
Root & Image\\\midrule
$\delta_1,\delta_2$ & $\alpha_7,\alpha_8$\\
$\delta_3$ & $-\alpha_1-2\alpha_2-3\alpha_3-2\alpha_4-\alpha_5-2\alpha_6$\\
$\delta_4$ & $\alpha_1$\\
$\delta_5$ & $\alpha_3+\alpha_4+\alpha_5+\alpha_6$\\
$\delta_6,\delta_7,\delta_8$ & $\alpha_2,\alpha_3,\alpha_4$\\\bottomrule
\end{tabular}
\end{center}
We also denote by $P=(P_{ij})$ the matrix of this map with respect to the ordered bases $(\delta_1,\ldots,\delta_8)$ of $L_{\rm Rie}$ and $(\alpha_1,\ldots,\alpha_8)$ of $L$, so that
\[
 P(\delta_j)=\sum_{i=1}^{8}P_{ij}\alpha_i
 \qquad(1\le j\le8).
\]
This matrix has integer entries and satisfies
\[
 \det P=-1,\qquad P^TGP=G_{\rm Rie},
\]
where $G$ is the Gram matrix in~\eqref{eq:grammatrix}. Thus $P$ is a lattice isometry. Here unimodularity refers to the matrix of $P$; both root lattices have Gram determinant $9$.

To obtain an isomorphism of lattice vertex operator algebras, we must also account for the cocycles associated with the two ordered bases. The cocycles on $L_{\rm Rie}$ and $L$ are $\eps_{\rm Rie}$ and $\eps$, respectively. These choices need not satisfy $\eps_{\rm Rie}(x,y)=\eps(Px,Py)$, so the assignment $e^x\mapsto e^{Px}$ requires a sign correction. For $x=\sum_{i=1}^8x_i\delta_i$, define
\begin{equation}\label{eq:Rlift}
 \Psi(e^x)=(-1)^{F(x)}e^{Px},\qquad
 F(x)=x_5+x_3(x_4+x_5+x_6+x_7+x_8)+x_5(x_6+x_7)\pmod2.
\end{equation}
The function $F$ satisfies
\[
 \eps_{\rm Rie}(x,y)(-1)^{F(x+y)}
 =(-1)^{F(x)+F(y)}\eps(Px,Py)
 \qquad(x,y\in L_{\rm Rie}),
\]
as follows by substitution in the two cocycle formulas. This identity ensures that $\Psi$ preserves the twisted group-algebra products. On the Heisenberg generators, set
\[
 \Psi\bigl(\delta_i(-1)\one\bigr)=(P\delta_i)(-1)\one
 \qquad(1\le i\le8).
\]
Together with the isometry property of $P$, these formulas give a vertex operator algebra isomorphism $\Psi:V_{L_{\rm Rie}}\longrightarrow V_L$ lifting $P$.

We now compare the displayed generators. The map $P$ satisfies
\[
 P\delta_3=\ell,\qquad Ph_{\rm Rie}=h,\qquad
 (P\zeta_1,\ldots,P\zeta_7)
 =(\alpha_7,\alpha_8,\alpha_1,\alpha_2,\alpha_3,\alpha_4,\alpha_5).
\]
For $x=\pm\delta_3$ and $x=\pm\zeta_i$, $1\le i\le7$, we have $F(x)\equiv0\pmod2$, and hence $\Psi(e^x)=e^{Px}$. It follows that
\[
 \Psi(H_\bullet^{\rm Rie})=H_\bullet,\qquad
 \Psi(E_\bullet^{\rm Rie})=E_\bullet,\qquad
 \Psi(F_\bullet^{\rm Rie})=F_\bullet
 \qquad(\bullet=L,R).
\]
In particular $\Psi(\omega_L^{\rm Rie})=\omega_L$ and $\Psi(\omega_R^{\rm Rie})=\omega_R$, while the lattice isometry gives $\Psi(\omega_{\rm Rie})=\omega$. For the three exponentials in $v$ one has
\[
 P\beta_1=c,\quad P\beta_2=d,\quad P\beta_3=a,
 \qquad
 \bigl((-1)^{F(\beta_1)},(-1)^{F(\beta_2)},(-1)^{F(\beta_3)}\bigr)=(-1,1,1).
\]
It follows term by term that
\[
 \Psi(v)=\frac{r-3\xi}{20}e^c
          +\frac{r+3\xi}{20}e^d+e^a=u.
\]

The four vectors $b_k$ can also be transported directly. If $z=\sum_i n_i\zeta_i$, substitution in~\eqref{eq:Rlift} gives
\[
 F(z)\equiv n_7(n_3+n_4+n_5+n_6)\pmod2.
\]
Moreover $\kappa_H(P\zeta_i)=m_i$. Consequently
\[
 \Psi\bigl(\mathsf S_{\rm Rie}(z)\bigr)
 =e^{Pz}+\kappa_H(Pz)e^{-Pz}.
\]
Together with the Heisenberg transformation, this proves that $B_k=\Psi(b_k)$ are exactly the vectors specified by~\eqref{eq:Btables} and Tables~\ref{tab:Bquadratic}--\ref{tab:Bexponential}. In particular the images of $e_R,f_R$ are completely determined by those tables. It remains to prove that the images of Riesen's two Virasoro generators belong to $U$. Put
\[
 e_*=\Psi(e_R)
 =\frac{2B_0+6B_1}{31}
   +\frac{2-r}{93}B_2+\frac{2-2r}{93}B_3,
 \qquad f_*=\Psi(f_R).
\]
Then $e_*+f_*=\omega_c$, where $\omega_c=\omega-\omega_L-\omega_R$. Thus $\Psi(U_{\rm Riesen})$ already contains the six affine generators, $u$, and $\omega$, and hence contains $U$. To prove the reverse inclusion, we shall recover $e_*$ and $f_*$ separately from Huang's generators. We do this by taking the component of $u_1\theta_H(u)$ annihilated by all six zero modes of the two $\slc$-triples.

We first check that $\theta_H(u)$ belongs to $U$. Since $\theta_H$ is antilinear and
\[
 \kappa_H(a)=4,\qquad \kappa_H(d)=\kappa_H(c)=40,
\]
its definition gives
\[
 \theta_H(u)=4e^{-a}+2(r-3\xi)e^{-d}+2(r+3\xi)e^{-c}.
\]
To obtain this vector from the generators of $U$, use $F_R=-\sum_{i\in I}h_i e^{-\alpha_i}$ and the recurrence
\[
 \frac{(F_R)_0^ku}{k!}
 =-\frac1k\sum_{i\in I}h_i(e^{-\alpha_i})_0
       \frac{(F_R)_0^{k-1}u}{(k-1)!}
 \qquad(1\le k\le13).
\]
Applying the lattice mode formula with the cocycle in~\eqref{eq:cocycle} to the three exponential terms separately gives
\[
 \frac{(F_R)_0^{13}}{13!}e^a=-34560\,e^{\ell-a},\quad
 \frac{(F_R)_0^{13}}{13!}e^d=-345600\,e^{\ell-c},\quad
 \frac{(F_R)_0^{13}}{13!}e^c=-345600\,e^{\ell-d}.
\]
Substituting these identities in~\eqref{eq:Hgenerator} gives
\[
 \frac{(F_R)_0^{13}u}{13!}
 =-8640\bigl(4e^{\ell-a}+2(r-3\xi)e^{\ell-d}
                         +2(r+3\xi)e^{\ell-c}\bigr).
\]
For $x=a,d,c$, the same formula gives $(F_L)_0e^{\ell-x}=e^{-x}$. Consequently,
\[
 \theta_H(u)=-\frac{(F_L)_0(F_R)_0^{13}u}{8640\cdot13!}\in U.
\]
Here $8640=\kappa_H(\ell)^{-1/2}=5\cdot8^2\cdot27$. Formula~\eqref{eq:extraction} for the coefficients also gives each step of the recurrence explicitly.

The vectors $e_*$ and $f_*$ have conformal weight two and are annihilated by all six zero modes of the two $\slc$-triples. To recover them from Huang's generators, we consider $u_1\theta_H(u)$, which has conformal weight two and Cartan weight zero for both $\slc$-triples. However, Cartan weight zero alone does not imply annihilation by $E_0$ and $F_0$. We therefore use Casimir operators to obtain its component annihilated by all six zero modes. For a zero-mode $\slc$-triple, put
\[
 \mathcal C(E,F,H)=H_0^2+2E_0F_0+2F_0E_0,
\]
and denote the Casimir operators for the left and right triples by $\mathcal C_L$ and $\mathcal C_R$, respectively. On the irreducible module $L(j)$ of highest weight $j$, the Casimir acts as $j(j+2)$, and the Cartan weights are $j,j-2,\ldots,-j$. Thus $L(j)$ contains weight zero exactly when $j$ is even.

A lattice monomial of conformal weight two has momentum $x$ with $x^2\le4$. Hence
\[
 |(\ell,x)|\le\sqrt8<3,\qquad |(h,x)|\le\sqrt{312}<18.
\]
The left and right highest weights in $(V_L)_2$ are therefore at most $2$ and $17$, respectively. On their zero-weight spaces, only the highest weights $0,2$ on the left and $0,2,4,\ldots,16$ on the right can contribute. Set
\begin{equation}\label{eq:casimirprojectors}
 \Pi_L=1-\frac{\mathcal C_L}{8},\qquad
 \Pi_R=\prod_{j=2,4,\ldots,16}\left(1-\frac{\mathcal C_R}{j(j+2)}\right).
\end{equation}
Each polynomial is $1$ on the trivial module and zero on every nontrivial module that meets the corresponding zero-weight space. Complete reducibility gives the projection property. The Casimir operators are self-adjoint for the Hermitian form in Section~4.4, so these projections are orthogonal. Since the two zero-mode actions commute, the projections also commute.

\begingroup
The projection property above is used only on the common zero-Cartan-weight subspace of $(V_L)_2$. The vectors $\omega_L$ and $\omega_R$ are also annihilated by all six zero modes of the two $\slc$-triples and are therefore retained by the projection, although they do not belong to the affine commutant. The calculation below gives their contributions as $6\omega_L+10\omega_R$. Subtracting these terms gives a linear combination of $e_*$ and $f_*$.
\par\endgroup
Accordingly, put
\begin{equation}\label{eq:cosetprojection}
 \widetilde u=\Pi_L\Pi_R\bigl(u_1\theta_H(u)\bigr),\qquad
 v_0=\widetilde u-6\omega_L-10\omega_R,
\end{equation}
and define
\begin{equation}\label{eq:recoveredvirasoro}
 e=\frac{35v_0-63\omega_c}{93},\qquad
 f=\omega_c-e.
\end{equation}
Here $\omega_c$ is the conformal vector of $C_U(A_{1,1}\otimes A_{1,39})$. All these vectors belong to $U$: we have proved $\theta_H(u)\in U$, and the projections are polynomials in zero modes of its affine generators. The next proposition identifies $e,f$ with $e_*,f_*$ and explains the coefficients in~\eqref{eq:recoveredvirasoro}.

\begin{proposition}\label{prop:RiesenEquality}
The vectors in \eqref{eq:recoveredvirasoro} satisfy $e=\Psi(e_R)$ and $f=\Psi(f_R)$. In particular,
\[
 \Psi(U_{\rm Riesen})=U.
\]
\end{proposition}
\begin{proof}
We first determine the components of $\widetilde u$ along $\omega_L,\omega_R,e_*,f_*$, and then prove that no further component remains. All inner products in this proof are taken with respect to $(\cdot,\cdot)_{\kappa_H}$. Put
\[
 a_0=\Pi_L\bigl(u_1\theta_H(u)\bigr),\qquad
 \widetilde u=\Pi_Ra_0.
\]
The four vectors $\omega_L,\omega_R,e_*,f_*$ are mutually orthogonal and are annihilated by the right affine zero modes, as follows from their coefficient formulas and the Gram matrix~\eqref{eq:invariantgram}. Hence, for each of these vectors $b$,
\[
 (\widetilde u,b)=(\Pi_Ra_0,b)=(a_0,\Pi_Rb)=(a_0,b).
\]
Substituting the expressions for $e_*,f_*$ in the $a_0$-column of~\eqref{eq:sourcepairings} gives the following values. The squared norms follow from~\eqref{eq:invariantgram}.
\[
\begin{array}{c|ccc}
 b&(b,b)&(\widetilde u,b)&(\widetilde u,b)/(b,b)\\\hline
 \omega_L&1/2&3&6\\
 \omega_R&117/82&585/41&10\\
 e_*&35/31&156/31&156/35\\
 f_*&1200/1271&2160/1271&9/5
\end{array}
\]
For example, the second coefficient is $(585/41)/(117/82)=10$. It follows that the orthogonal projection of $\widetilde u$ onto the span of these four vectors is
\[
 T_0=6\omega_L+10\omega_R+\frac{156}{35}e_*+\frac95f_*.
\]
To prove $\widetilde u=T_0$, we must also rule out a component orthogonal to all four vectors. Since $\Pi_R$ is an orthogonal projection, the last entry in the $a_0$-column of~\eqref{eq:projectionmoments} gives
\[
 (\widetilde u,\widetilde u)=(a_0,\Pi_Ra_0)=\frac{6516}{35}.
\]
On the other hand, the preceding table and orthogonality give
\[
 \begin{aligned}
 (T_0,T_0)
 &=6^2\frac12+10^2\frac{117}{82}
   +\left(\frac{156}{35}\right)^2\frac{35}{31}
   +\left(\frac95\right)^2\frac{1200}{1271}\\
 &=\frac{6516}{35}.
 \end{aligned}
\]
Thus $\|\widetilde u-T_0\|^2=\|\widetilde u\|^2-\|T_0\|^2=0$, and positive definiteness proves
\[
 \widetilde u=6\omega_L+10\omega_R+\frac{156}{35}e_*+\frac95f_*.
\]
This argument does not assume that $\widetilde u$ lies in the span of the four vectors. The lattice mode recursions and all scalar data used above are given in Appendix~\ref{app:checks}.

Subtracting the two Sugawara terms and using $e_*+f_*=\omega_c$, we obtain
\[
 v_0=\frac{156}{35}e_*+\frac95f_*
     =\frac{93}{35}e_*+\frac{63}{35}\omega_c.
\]
Solving for $e_*$ gives precisely~\eqref{eq:recoveredvirasoro}. Hence $e=e_*$ and $f=f_*$. In particular,
\begin{equation}\label{eq:riesenaxis}
 e=\frac2{31}\Psi(b_0)+\frac6{31}\Psi(b_1)
       +\frac{2-r}{93}\Psi(b_2)+\frac{2-2r}{93}\Psi(b_3).
\end{equation}

Every displayed generator of $\Psi(U_{\rm Riesen})$ now belongs to $U$: the six affine generators and $u$ by the explicit correspondence, and $e_*,f_*$ by their construction above. Thus $\Psi(U_{\rm Riesen})\subseteq U$. Conversely, the image contains the six affine generators, $u$, and
\[
 \omega=\omega_L+\omega_R+e+f.
\]
These generate $U$, so $U\subseteq\Psi(U_{\rm Riesen})$. The two inclusions prove the asserted equality, and restriction of the lattice vertex operator algebra isomorphism $\Psi$ gives $U_{\rm Riesen}\cong U$.
\end{proof}

For completeness, the following products give an independent check of the recovered Virasoro vectors. The lattice calculation gives, for each affine generator $a$,
\[
 a_ne=a_nf=0\quad(n=0,1,2),
\]
and
\begin{align}\label{eq:efproducts}
 e_0e&=\partial e,& e_1e&=2e,&e_2e&=0,&e_3e&=\frac{35}{31}\one,\notag\\
 f_0f&=\partial f,& f_1f&=2f,&f_2f&=0,&f_3f&=\frac{1200}{1271}\one.
\end{align}
All mixed products $e_nf,f_ne$ vanish for $n\ge0$. Also,
\[
 e_1u=\frac{13}{31}u,\qquad f_1u=\frac{180}{1271}u,\qquad
 v_0=\frac{156}{35}e+\frac95f.
\]
Thus $e,f$ are commuting Virasoro vectors of central charges $70/31$ and $2400/1271$. Appendix~\ref{app:invariants} gives the coefficient tables and the positive-definite Gram-matrix calculations for these identities.

\section{The lattice neighbor and the two involutions}\label{sec:neighbor}
We construct two even lattices isometric to $E_6\oplus A_2$ with a common sublattice of index two in each. The associated lattice-character involutions have isomorphic fixed point lattice vertex operator algebras. We also verify that these involutions preserve the subalgebras generated by Huang's and Jiao's vectors, respectively.

\begingroup
We now use two coordinate copies of the lattice in Section~\ref{sec:lattice}. We retain the notation $L$ for Huang's original copy and write its ordered root basis as $\alpha_1^H,\ldots,\alpha_8^H$. We denote Jiao's copy by $L_J$, with ordered root basis $\alpha_1^J,\ldots,\alpha_8^J$. Thus
\[
 L=\bigoplus_{i=1}^8\Z\alpha_i^H,\qquad
 L_J=\bigoplus_{i=1}^8\Z\alpha_i^J,\qquad
 (\alpha_i^H,\alpha_j^H)=(\alpha_i^J,\alpha_j^J)=G_{ij}.
\]
The superscripts distinguish the two copies, which use the same root numbering and Gram matrix. The formulas for Huang's generators are read in the $H$ copy, and those for Jiao's generators in the $J$ copy, with all coefficients unchanged. In particular, $\ell,h\in L$, whereas $p,b\in L_J$. In the Jiao coordinate space, define
\par\endgroup
\begin{equation}\label{eq:neighborcharges}
 \lambda=\frac{p+b}{2},\qquad s=\frac{p+b}{4},\qquad t=\frac{13p-3b}{4}.
\end{equation}
Then
\begin{equation}\label{eq:norms}
 \lambda^2=8,\quad s^2=2,\quad t^2=78,\quad(s,t)=0,
 \quad p=\frac{3s+t}{4},\quad b=\frac{13s-t}{4}.
\end{equation}
In particular $(s,p)=3/2$. Put
\begin{equation}
 L_0=\{x\in L_J:(\lambda,x)\in2\Z\},\qquad L_H=L_0+\Z s.
\end{equation}
\begingroup
The lattice $L_H$ lies in $L_J\otimes_{\Z}\Q$. We construct a lattice isometry $T:L\longrightarrow L_H$ in the proof of Proposition~\ref{prop:lattice} below.
\par\endgroup
Since $G\lambda=(1,1,-1,1,1,-1,1,1)^T$, $L_0$ consists of the integer coordinate vectors with even coordinate sum. It has index two in $L_J$, and $p$ represents the odd coset.

\begin{proposition}\label{prop:lattice}
$L_H$ is an even lattice isometric to $E_6\oplus A_2$ and hence isometric to $L_J$. The decompositions into $L_0$-cosets are
\[
 L_J=L_0\sqcup(p+L_0),\qquad L_H=L_0\sqcup(s+L_0).
\]
Here $\sqcup$ denotes disjoint union.
\end{proposition}
\begin{proof}
We have $2s=\lambda\in L_0$, $s\notin L_J$, $(s,L_0)\subset\Z$ and $s^2=2$. This proves integrality and evenness. Both $L_J$ and $L_H$ have determinant nine. To identify the latter without relying on determinant alone, let $e_i=\alpha_i^J$ and use the ordered basis
\begin{equation}\label{eq:BHbasis}
 \mathcal B_H=(2e_1,e_2-e_1,e_3+e_1,e_4-e_1,e_5-e_1,s,e_7-e_1,e_8-e_1).
\end{equation}
In this basis take the eight column vectors
\begin{equation}\label{eq:rtable}
\begin{array}{ll}
 r_1=(2,2,1,1,1,-1,0,0)^T,&r_2=(2,1,1,1,0,-1,1,1)^T,\\
 r_3=(-6,-3,-2,-3,-2,3,-2,-2)^T,&r_4=(2,1,1,1,1,-1,1,1)^T,\\
 r_5=(2,1,1,2,1,-1,0,0)^T,&r_6=(3,1,0,1,1,-2,1,1)^T,\\
 r_7=(2,1,1,1,1,-1,1,0)^T,&r_8=(2,1,1,1,1,-1,0,1)^T.
\end{array}
\end{equation}
Their Gram matrix is $G$, and their determinant relative to $\mathcal B_H$ has absolute value one. They are therefore a basis, not merely a root sublattice.
\begingroup
Consequently the map
\[
 T:L\longrightarrow L_H,\qquad T(\alpha_i^H)=r_i
 \quad(1\le i\le8)
\]
is a lattice isometry. Here the $r_i$ form a root basis of $L_H$ in the Jiao coordinate space, whereas the $\alpha_i^J$ form the chosen basis of $L_J$. The isometry satisfies
\par\endgroup
\[
 T\ell=s,\qquad Th=t.
\]
All these assertions follow by multiplying the displayed rational matrices.
\end{proof}

On $V_{L_H}$ and $V_{L_J}$ respectively define
\begin{equation}\label{eq:involutions}
 \tau=\exp(2\pi i p(0)),\qquad \sigma=\exp(\pi i\lambda(0)).
\end{equation}
The first is $+1$ on $L_0$ and $-1$ on $s+L_0$; the second is $+1$ on $L_0$ and $-1$ on $p+L_0$. Hence
\begin{equation}\label{eq:ambientfixed}
 (V_{L_H})^\tau=V_{L_0}=(V_{L_J})^\sigma.
\end{equation}
A cocycle-compatible identification is understood here. On the Huang coordinate copy,
\[
 \tau=\exp\!\left(\frac{\pi i}{2}(3H_L+H_R)_{0}\right).
\]
It fixes $H_L,H_R,\omega,u$ and negates the four root generators. On the Jiao generators $\sigma$ fixes the two Cartans, the conformal vector and both weight-two generators, and negates the four root generators. Thus both generated subalgebras are invariant.

\begin{remark}
Writing both presentations in the same root convention does not make their neighbor comparison the identity embedding in $V_L$. The comparison uses the rational isometry $T$ and its lift on the common sublattice. In particular, \eqref{eq:ambientfixed} does not imply $U^\tau=V^\sigma$ for the smaller generated subalgebras.
\end{remark}

\section{The lattice isomorphism and the correspondence of generators}\label{sec:lift}
We construct the lattice vertex operator algebra isomorphism used in the common fixed point comparison. The calculation determines the cocycle and character factors, explains the required interchange of the two $A_2$ roots, and gives the explicit correspondence of generators. We also verify compatibility with the unitary structures.

\subsection{Lifts of the lattice isomorphism}
\leavevmode
\begingroup
Let $R=(R_{ji})$ be the matrix determined by
\[
 T(\alpha_i^H)=r_i=\sum_{j=1}^{8}R_{ji}\alpha_j^J
 \qquad(1\le i\le8).
\]
Thus $[T(y)]_{\alpha^J}=R[y]_{\alpha^H}$, where the brackets denote coordinate columns in the indicated ordered bases. The columns of $R$ are obtained by expressing the vectors $r_i$ in \eqref{eq:rtable} in the $\alpha^J$ coordinates. Let $\eps_H$ and $\eps_J$ denote the cocycles defined by \eqref{eq:cocycle} on $L$ and $L_J$, respectively, using their ordered root bases. We seek a lift of $R:R^{-1}L_0\longrightarrow L_0$ to an isomorphism of lattice vertex operator algebras. For $x,y\in L_0$, the cocycle factors in the source and target products are $\eps_H(R^{-1}x,R^{-1}y)$ and $\eps_J(x,y)$, respectively. The scalar factors in the lift must satisfy \eqref{eq:liftidentity} below to preserve the twisted group-algebra products.

To specify one such choice, replace $s$ by $2s$ in the ordered basis $\mathcal B_H$ of \eqref{eq:BHbasis}. This gives a basis of $L_0$: the other seven vectors lie in $L_0$, and $s+L_0$ generates $L_H/L_0$ of order two. Let $\mathsf B$ be the matrix whose columns are this basis in the $\alpha^J$ coordinates, and let $\mathbf e_i$ be the $i$th standard coordinate column of $\Z^8$, for $1\le i\le8$. For $x=\mathsf Bn\in L_0$, where $n=(n_1,\ldots,n_8)^T\in\Z^8$, set
\begin{equation}\label{eq:gamma}
 \gamma(x)=\prod_{i>j}
 \left[\eps_H(R^{-1}\mathsf B\mathbf e_i,R^{-1}\mathsf B\mathbf e_j)
              \eps_J(\mathsf B\mathbf e_i,\mathsf B\mathbf e_j)\right]^{n_i n_j}.
\end{equation}
Since $R$ is an isometry, the cocycles $\eps_H(R^{-1}x,R^{-1}y)$ and $\eps_J(x,y)$ on $L_0$ have equal diagonal values and equal commutators. It follows that
\begin{equation}\label{eq:liftidentity}
 \gamma(x)\gamma(y)\eps_J(x,y)
 =\eps_H(R^{-1}x,R^{-1}y)\gamma(x+y).
\end{equation}
This identity shows that
\[
 \phi_0(e^{R^{-1}x})=\gamma(x)e^x
\]
preserves the twisted group-algebra products. Together with the Heisenberg map induced by $R$, it defines an isomorphism $\phi_0:V_{R^{-1}L_0}\longrightarrow V_{L_0}$.

Any other lift of the same lattice map is
\begin{equation}\label{eq:alllifts}
 \phi_\chi(e^{R^{-1}x})=\gamma(x)\chi(x)e^x,
 \qquad\chi\in\operatorname{Hom}(L_0,\C^\times).
\end{equation}
The lattice character $\chi$ is uniquely determined by its nonzero values on the eight columns of $\mathsf B$, since
\[
 \chi(\mathsf Bn)=\prod_{i=1}^8\chi(\mathsf B\mathbf e_i)^{n_i}.
\]
These values describe the remaining freedom when matching the coefficients of the generators in the next subsection.
\par\endgroup

\subsection{The choice of lattice isomorphism}
\leavevmode
\begingroup
We now ask whether a lift of $R$ can identify certain weight-two vectors in Huang's realization with those in Jiao's realization. These correspondences will be used to compare the fixed-point algebras. We first compare the lattice vectors occurring in their exponential expressions. For a finite sum $v=\sum_x c_xe^x$, its support is
\[
 \operatorname{supp}(v)=\{x:c_x\ne0\}.
\]
Direct calculation gives
\[
 \begin{aligned}
 R\bigl(\operatorname{supp}((E_L)_{-1}E_R)\bigr)
   &=\operatorname{supp}(u_{2,6}),\\
 R\bigl(\operatorname{supp}((E_L)_{-1}F_R)\bigr)
   &=\operatorname{supp}(u_{1,7}).
 \end{aligned}
\]
In each case, both sides consist of seven lattice vectors. The three lattice vectors $a,d,c$ occurring in $u$ satisfy
\par\endgroup
\begin{align}
 Ra&=\rho_1+\rho_3,& Rd&=\rho_2+\rho_3,& Rc&=\rho_1+\rho_2.
\end{align}
\begingroup
Since the $\rho_i$ are pairwise orthogonal, the support of $(P_E)_{-1}P_E$ consists of these three lattice vectors. The following proposition shows that equality of these supports does not determine the required correspondence of generators.
\par\endgroup

\begin{proposition}[A restriction on the lift of $R$]\label{prop:obstruction}
For the lattice map $R$ above and the cocycles and generators in Section~\ref{sec:lattice}, there is no lift $\phi_\chi$ and no nonzero scalars $A,D$ such that both
\begin{equation}\label{eq:impossible}
 \phi_\chi((E_L)_{-1}E_R)=A u_{2,6},
 \qquad \phi_\chi(u)=D(P_E)_{-1}P_E
\end{equation}
hold. This remains true when arbitrary complex lattice-character values are allowed.
\end{proposition}
\begin{proof}
Write $m_i=R(\ell+\alpha_i^H)$ for $i\in I$, and $y_0=Ra$, $y_1=Rd$. Their coordinates give
\begin{equation}\label{eq:momentumrelation}
 m_4+m_5+y_0=m_7+m_8+y_1.
\end{equation}
Let $d_i$ denote the coefficient of $e^{m_i}$ in $u_{2,6}$, and put
\[
 t_i=\frac{d_i}{\eps_H(\ell,\alpha_i^H)\gamma(m_i)}.
\]
For the four indices needed in \eqref{eq:momentumrelation}, direct substitution gives
\begin{equation}\label{eq:ti}
 \begin{array}{c|cccc}
 i&4&5&7&8\\\hline
 t_i&(r+\xi )/6&-(1-r\xi )/6&-(1+\xi )/3&(1-\xi )/6.
 \end{array}
\end{equation}
The first identity in \eqref{eq:impossible} would imply $\chi(m_i)=At_i$. Applying the character to \eqref{eq:momentumrelation} therefore forces
\begin{equation}\label{eq:forcedratio}
 \frac{\chi(y_1)}{\chi(y_0)}
 =\frac{t_4t_5}{t_7t_8}=\frac{r-3\xi }{2}.
\end{equation}
The factor $A$ cancels.

On the other hand the coefficients of $(P_E)_{-1}P_E$ at $y_0,y_1$ are both $-2$, while $\gamma(y_0)=1$ and $\gamma(y_1)=-1$. The second identity in \eqref{eq:impossible}, using the coefficients $1,(r+3\xi )/20$ in $u$, would force
\begin{equation}\label{eq:requestedratio}
 \frac{\chi(y_1)}{\chi(y_0)}
 =-\frac{20}{r+3\xi }=-\frac{r-3\xi }{2},
\end{equation}
since $(r+3\xi )(r-3\xi )=40$. Equations \eqref{eq:forcedratio} and \eqref{eq:requestedratio} contradict one another.
\end{proof}

We next choose a different isometry for which the required identities hold.

\subsection{The isomorphism used in the comparison}
Let $D$ fix $\alpha_1,\ldots,\alpha_6$ and interchange $\alpha_7,\alpha_8$. Set
\begin{equation}\label{eq:correctedmap}
 R_* = RD,\qquad R_*\alpha_i=r_i\ (1\le i\le6),\quad
 R_*\alpha_7=r_8,\quad R_*\alpha_8=r_7.
\end{equation}
The symbol $R_*$ denotes the composition $R\circ D$; the star distinguishes this map from $R$ and does not denote an adjoint. Thus $R_*$ is obtained by interchanging the seventh and eighth columns of $R$. It satisfies $R_*\ell=s$ and $R_*h=t$. A standard lift of $D$ is
\[
 \widehat D(e^x)=(-1)^{x_7x_8}e^{Dx}.
\]
It preserves the two $\slc$-triples but sends
\[
 u\longmapsto-e^a+\frac{r+3\xi }{20}e^d+\frac{r-3\xi }{20}e^c.
\]
The lift therefore changes the relative sign of the first summand of $u$. We use $R_*$ in the comparison below.

Use the same ordered basis $\mathsf B$ of $L_0$. For $x=\mathsf Bn$ define $\gamma_*(x)=(-1)^{q_*(n)}$, where
\begin{align}\label{eq:qstar}
 q_*(n)={}&n_2n_1+n_5(n_1+n_2+n_3+n_4)+n_6(n_2+n_5)\notag\\
       &+n_7(n_2+n_3+n_4+n_6)+n_8n_1\pmod2.
\end{align}
This is exactly formula~\eqref{eq:gamma} with $R$ replaced by $R_*$. Define a character $\chi_*$ by its values on the eight columns $\mathbf b_i=\mathsf B\mathbf e_i$ of $\mathsf B$:
\begin{equation}\label{eq:characterstar}
\begin{array}{c|l}
 i&\chi_*(\mathbf b_i)\\\hline
1&1/720\\
2&(r-7\xi )/8\\
3&(-537r+1377\xi )/128\\
4&(-r+3\xi )/2\\
5&-(5+r\xi )/16\\
6&-(7+5r\xi )/2388787200\\
7&(r-3+(1+r)\xi )/4\\
8&(-r-3+(r-1)\xi )/8
\end{array}
\end{equation}
Then the oscillator map induced by $R_*$ and
\begin{equation}\label{eq:phistar}
 \Phi(e^{R_*^{-1}x})=\gamma_*(x)\chi_*(x)e^x
\end{equation}
define the required common-lattice VOA isomorphism. These values specify one choice of normalization for the lift.

\begin{proposition}[The correspondence of generators]\label{prop:dictionary}
Put
\[
 c_{26}=\frac{7+5r\xi }{92160},\qquad c_{17}=\frac{\xi }{8640},\qquad
 c_P=\frac{32765+71r\xi }{32768}.
\]
The lift~\eqref{eq:phistar} satisfies
\begin{align}\label{eq:dictionary}
 \Phi((E_L)_{-1}E_R)&=c_{26}u_{2,6},&
 \Phi((E_L)_{-1}F_R)&=c_{17}u_{1,7},\notag\\
 \Phi(u)&=c_P\nn{P_E^2},&
 \Phi(\theta_Hu)&=c_P^{-1}\nn{P_F^2}.
\end{align}
Moreover $\Phi\theta_H=\theta_J\Phi$ on the common lattice VOA.
\end{proposition}
\begin{proof}
The cocycle identity~\eqref{eq:liftidentity} for $\gamma_*$ is a bilinear identity, verified on the 64 ordered basis pairs of $\mathsf B$. Substitution of~\eqref{eq:characterstar} into the seven, seven and three exponential coefficients gives the first three identities in~\eqref{eq:dictionary}. For completeness, the identity $\Phi\theta_H=\theta_J\Phi$ reduces on the same basis to
\[
 |\chi_*(x)|^2=\frac{\kappa_H(R_*^{-1}x)}{\kappa_J(x)},
\]
where $\kappa_H(x)=5^{x_1}8^{x_2}9^{x_3}8^{x_4}5^{x_5}2^{x_7+x_8}$ and
$\kappa_J(x)=2^{-x_1+x_2-x_4+x_5-x_7+x_8}$. Both sides are positive characters, so the eight checks prove the identity on all of $L_0$. The fourth identity follows by applying $\theta_J$ to the third and using $\Phi\theta_H=\theta_J\Phi$, $|c_P|=1$ and $\theta_J(P_E)=-P_F$. It also follows by direct substitution in the three coefficients.
\end{proof}

The stability of $U$ under $\theta_H$ and of $V$ under $\theta_J$ can be checked directly. In $U$ one has
\begin{equation}\label{eq:huangpct}
 (F_L)_{0}(F_R)_{0}^{13}u=-8640\cdot13!\,\theta_Hu.
\end{equation}
In $V$ the corresponding formulas are
\[
 (P_F)_{0}^2(B_F)_{0}^6u_{2,6}=1440\theta_Ju_{2,6},\qquad
 (P_F)_{0}(B_F)_{0}^7u_{1,7}=5040\theta_Ju_{1,7}.
\]
These identities follow from the lattice vertex-operator formula and, together with the relations among the weight-one generators, show that $\theta_H(U)=U$ and $\theta_J(V)=V$.

\begin{corollary}\label{cor:unitary}
The vertex operator algebras $U$ and $V$ are simple, unitary and self-contragredient, and are of CFT type.
\end{corollary}
\begin{proof}
Both algebras contain the lattice conformal vector and are stable under the corresponding antilinear involution. Restriction of the positive definite lattice form therefore makes them unitary. Their weight-zero spaces are $\C\one$. Complete reducibility of a unitary module, applied to the adjoint module, implies simplicity in CFT type \cite{DLUnitary}. The invariant Hermitian form and the antilinear involution give a nondegenerate invariant bilinear form, so both algebras are self-contragredient.
\end{proof}

The isomorphism $\Phi$ on the common lattice vertex operator algebra extends to an isomorphism from the Huang lattice realization onto $V_{L_H}$, after a compatible cocycle on $L_H$ is fixed. To see this, first choose any lift of $R_*$. The ratio of its restriction and $\Phi$ is a character of the index-two sublattice. Such a character extends to the full lattice by choosing a square root on the odd coset. The norm identity in Proposition~\ref{prop:dictionary} allows this extension to preserve the unitary structures. The two square-root choices differ by $\tau$ and do not change the fixed point comparison.

\section{Parafermions and the common fixed-point theorem}\label{sec:fixed}
We give explicit generators for the parafermion vertex operator algebras $K(\slc,3)$ and $K(\slc,13)$ in Jiao's realization. A lattice-shift lemma and a finite set of mode identities then prove that Huang's and Jiao's algebras have isomorphic fixed point subalgebras under the indicated involutions. This isomorphism also gives embeddings of the two parafermion algebras in Huang's fixed point subalgebra.

\subsection{Generators of the parafermion vertex operator algebras}
For any displayed level-$k$ $\slc$-triple $(E,F,H)$ define
\begin{equation}\label{eq:eta}
 \eta_k=\Omega_k(E,F,H)-\frac1{4k}H_{-1}H
 =\frac{-kH_{-2}\one-H_{-1}H+2kE_{-1}F}{2k(k+2)}.
\end{equation}
Here $H_{-2}\one$ denotes the derivative of the weight-one vector. A weight-three generator can be chosen as
\begin{align}\label{eq:W3}
 W_k^3={}&k^2H_{-3}\one+3kH_{-2}H+2H_{-1}^3\one
          -6kH_{-1}E_{-1}F\notag\\
 &+3k^2E_{-2}F-3k^2E_{-1}F_{-2}\one.
\end{align}
These are the standard affine formulas of \cite{DLLWY}. For $k\ge3$, $\eta_k,W_k^3$ generate the simple parafermion algebra $K(\slc,k)$; generation here is again vertex generation. Substitution of \eqref{eq:P}--\eqref{eq:BH} gives explicit lattice vectors for $k=3,13$. Their central charges are $4/5$ and $8/5$. The regularity of these standard positive-integral-level parafermion algebras follows from \cite{ALY}.

The Jiao involution has
\[
 \sigma=\exp\!\left(\frac{\pi i}{2}(P_H+B_H)_{0}\right).
\]
Each parafermion has zero charge under its own Cartan, and the two affine factors commute. Consequently
\begin{equation}\label{eq:parainfixed}
 K(\slc,3)\otimes K(\slc,13)\subset V^\sigma.
\end{equation}
This inclusion does not use rationality of $V$ or the fixed-point isomorphism proved below.

In Huang's construction the desired level-$3$ and level-$13$ parafermions are not the Cartan commutants of its level-$1$ and level-$39$ affine factors. Theorem~\ref{thm:fixed} below proves that their images under $\Phi^{-1}$ belong to $U^\tau$. Thus transporting \eqref{eq:eta} and \eqref{eq:W3} through the explicit lift~\eqref{eq:phistar} gives actual embedded parafermions in Huang's algebra. The identity $\Phi\theta_H=\theta_J\Phi$ also identifies the antilinear involutions on these subalgebras. Their central charges alone would not have established this membership.

\subsection{A lattice-shift lemma}
In this subsection let $M$ be a positive definite even lattice, and let $u\in M$ be nonzero; this local notation is independent of the weight-two generator used elsewhere. Fix a bimultiplicative lattice cocycle $\eps$ on $M$, with $\eps(x,y)\eps(y,x)^{-1}=(-1)^{(x,y)}$. For $i\in\Z$ and a lattice monomial $Pe^x$, where $P$ is a polynomial in the operators $h(-n)$ with $h\in\C\otimes_\Z M$ and $n>0$, set
\[
 S_{iu}(Pe^x)=\eps(iu,x)Pe^{x+iu}.
\]
This is a linear momentum-shift operator, not a VOA automorphism and not a weight-preserving map.

\begin{lemma}\label{lem:shift}
Let $A\subset V_M$ be a vertex subalgebra containing $V_{\Z u}$, and let $g$ be a lattice-character involution preserving $A$ and satisfying $g(e^u)=-e^u$. Suppose that $W$ is a vertex subalgebra satisfying
\[
 V_{2\Z u}\subset W\subset A^g
\]
and every chosen generator of $A$ lies in $W+S_uW$. Then
\[
 A=W\oplus S_uW,\qquad A^g=W.
\]
\end{lemma}
\begin{proof}
\leavevmode
\begingroup
We first show that every vector in $W$ is a linear combination of vectors obtained by applying finite products of the operators $u(-n)$, $n>0$, to vectors $a\in W$ satisfying
\[
 u(0)a=qa,\qquad u(n)a=0\quad(n>0),
\]
where $q\in\Z$.

Since $u(-1)\one\in V_{2\Z u}\subset W$, every mode $u(n)$ preserves $W$. On $Pe^x$, the operator $u(0)$ acts as multiplication by $(u,x)$. Each vector has only finitely many $u(0)$-eigencomponents, so polynomial interpolation in $u(0)$ shows that each component of a vector in $W$ also belongs to $W$.

For $n>0$, put
\[
 Q_n=\sum_{j\ge0}
       \frac{(-1)^j u(-n)^j u(n)^j}{j!(n u^2)^j}.
\]
The sum is finite on each vector because $u(n)$ is locally nilpotent. For $a\in V_M$ with $u(n)a=0$ and $m\ge0$, the relation $[u(n),u(-n)]=n u^2$ and the binomial theorem give
\[
 Q_n\bigl(u(-n)^m a\bigr)
 =\left(\sum_{j=0}^{m}(-1)^j\binom mj\right)u(-n)^m a
 =\begin{cases}
   a,&m=0,\\
   0,&m>0.
  \end{cases}
\]
Thus $Q_n$ retains precisely the terms with no $u(-n)$ factor. Now take $v\in W$ with $u(0)v=qv$. The Heisenberg decomposition of $V_M$ gives a finite expansion
\[
 v=\sum_{m=0}^{N}u(-n)^m a_m,\qquad a_m\in V_M,\quad u(n)a_m=0.
\]
The vector coefficients are given by
\[
 a_m=\frac{Q_nu(n)^m v}{m!(n u^2)^m}.
\]
All operators on the right preserve $W$ and commute with $u(0)$, so $a_m\in W$ and $u(0)a_m=qa_m$. Repeating this for the finitely many indices $n>0$ for which $u(-n)$ occurs in $v$ proves the asserted decomposition.
\par\endgroup

For an integer $i$ and such a vector $a\in W$ with $u(0)$-eigenvalue $q$, the lattice formula gives
\[
 Y(e^{iu},z)a=z^{iq}E^-(iu,z)S_{iu}a,
 \qquad E^-(iu,z)=\exp\!\sum_{n>0}\frac{iu(-n)}n z^n.
\]
\begingroup
Since $E^-(iu,z)$ has constant term $1$ and only nonnegative powers of $z$, the coefficient of the lowest power $z^{iq}$ in $Y(e^{iu},z)a$ is $S_{iu}a$. For $i=\pm2$, we have $e^{iu}\in V_{2\Z u}\subset W$, and hence
\[
 S_{2u}a=(e^{2u})_{-2q-1}a\in W,\qquad
 S_{-2u}a=(e^{-2u})_{2q-1}a\in W.
\]
The operators $S_{2u}$ and $S_{-2u}$ commute with every $u(-n)$, $n>0$, so the decomposition proved above gives $S_{2u}W\subset W$ and $S_{-2u}W\subset W$. Moreover,
\[
 S_{2u}S_{-2u}=\eps(2u,-2u)\operatorname{id}.
\]
Since this scalar is nonzero, $W\subset S_{2u}W$, and therefore $S_{2u}W=W$. For $i=\pm1$, the same argument, using $e^{\pm u}\in A$, gives $S_{\pm u}W\subset A$.
\par\endgroup

Now let $a,b\in W$ satisfy $u(0)a=qa$, $u(0)b=rb$ and $u(n)a=u(n)b=0$ for $n>0$. For integers $i,j$, the lattice formula gives
\[
 Y(S_{iu}a,z)S_{ju}b
 =(-1)^{jq}\eps(iu,ju)z^{ir+jq+ij u^2}
   E^-(iu,z)S_{(i+j)u}Y(a,z)b.
\]
To verify the scalar factor, write $a=Pe^x$ and $b=Qe^y$, with $P,Q$ involving only Heisenberg modes orthogonal to $u$. Bimultiplicativity and the commutator of $\eps$ give
\[
 \frac{\eps(iu,x)\eps(ju,y)\eps(x+iu,y+ju)}
      {\eps(x,y)\eps((i+j)u,x+y)}
 =(-1)^{jq}\eps(iu,ju).
\]
The added momentum pairing is $ir+jq+ij u^2$, and the perpendicular Heisenberg modes commute with the shifts. Linearity gives the formula for arbitrary $a,b$ satisfying the stated conditions.

Also $S_{iu}S_{ju}=\eps(iu,ju)S_{(i+j)u}$. Since $S_{2u}W=W$, the space $S_{ku}W$ equals $W$ for even $k$ and $S_uW$ for odd $k$. Every coefficient of $Y(a,z)b$ lies in $W$, and the coefficients of $E^-(iu,z)$ are polynomials in the operators $u(-n)$, $n>0$. The displayed product therefore has all coefficients in $W+S_uW$. This space is preserved by every $u(n)$, so the Heisenberg commutator and normal-ordering formulas extend this conclusion to oscillator descendants. Thus $W+S_uW$ is a vertex subalgebra containing the generators of $A$, and equals $A$. The two summands have opposite $g$-parity.
\end{proof}

\subsection{Generators of the fixed point subalgebras}
Use the lift $\Phi$ on $L_0$ and transport Huang's presentation into $V_{L_H}$. Denote the transported algebra by $\widetilde U$, and suppress tildes on its generators. Inside its even part let
\begin{equation}\label{eq:Wdef}
 W=\langle\omega,H_L,H_R,e^{2s},e^{-2s},u,
                 (E_L)_{-1}E_R,(E_L)_{-1}F_R\rangle_{\rm VA}.
\end{equation}
The exponentials $e^{\pm2s}$ belong to the rank-one lattice algebra generated by $E_L,F_L,H_L$. All vectors in \eqref{eq:Wdef} are even. The right root generators are recovered by shifting the last two vectors by $-s$ and taking their parts with no Heisenberg factors. Lemma~\ref{lem:shift} therefore gives, without a comparison to $V$,
\begin{equation}\label{eq:Ueven}
 \widetilde U=W\oplus S_sW,\qquad\widetilde U^\tau=W.
\end{equation}

\begin{proposition}[Finite sufficient conditions for the comparison]\label{prop:comparison}
Suppose the chosen lift satisfies the following explicit membership conditions:
\begin{enumerate}
\item Every vector in \eqref{eq:Wdef} belongs to $V^\sigma$.
\item The vectors $u_{2,6},u_{1,7},e^{\pm2p},\nn{P_E^2},\nn{P_F^2}$ and the four vectors $\nn{P_\varepsilon B_\delta}$, $\varepsilon,\delta\in\{+,-\}$, belong to $W$, where $P_+=P_E$, $P_-=P_F$, $B_+=B_E$ and $B_-=B_F$.
\end{enumerate}
Then $\widetilde U^\tau=W=V^\sigma$.
\end{proposition}
\begin{proof}
The first condition gives $W\subset V^\sigma$. The Cartans $p,b$ are linear combinations of $s,t$, so $P_H,B_H\in W$. Since the three $\rho_i$ are pairwise orthogonal, direct normal multiplication gives
\begin{equation}\label{eq:cubic}
 \nn{P_E^3}=-6e^p,\qquad \nn{P_F^3}=6e^{-p}.
\end{equation}
Thus $V_{\Z p}\subset V$, and the second condition gives $V_{2\Z p}\subset W$.
The shifts $S_{-p}P_E$ and $S_pP_F$ are nonzero scalar multiples of $\nn{P_F^2}$ and $\nn{P_E^2}$, respectively. Since the two affine algebras commute, their positive root generators have regular OPE with one another. Consequently
\[
 (\nn{P_E^2})_{-1}(\nn{P_EB_\delta})
       =\nn{P_E^3B_\delta}
       =-6\nn{e^pB_\delta}.
\]
The lattice vectors occurring in $B_\delta$ are orthogonal to $p$, so the last expression is a nonzero scalar multiple of $S_pB_\delta$. The negative formulas are analogous. Every root generator of $V$ therefore lies in $W+S_pW$, and its remaining defining generators lie in $W$. Lemma~\ref{lem:shift} yields $V=W\oplus S_pW$ and $V^\sigma=W$.
\end{proof}

\subsection{The common fixed point subalgebra}
\begin{theorem}[Common fixed point subalgebra]\label{thm:fixed}
For the lift $\Phi$, the sufficient conditions in Proposition~\ref{prop:comparison} hold. Consequently
\[
 \Phi(U^\tau)=W=V^\sigma,\qquad
 \widetilde U=W\oplus S_sW,\quad V=W\oplus S_pW.
\]
The assertion uses only the defining vectors, exact mode identities and Lemma~\ref{lem:shift}; it does not assume rationality of $U$ and $V$ or equality of characters.
\end{theorem}
\begin{proof}
By Proposition~\ref{prop:dictionary}, the three non-Cartan weight-two generators of $W$ belong to $V^\sigma$. Direct multiplication gives
\begin{equation}\label{eq:e2s}
 (u_{2,6})_{-1}u_{1,7}=-\frac{13\xi }{3}e^{2s}.
\end{equation}
Since $\theta_J(V)=V$ and $\theta_J(e^{2s})=\kappa_J(2s)e^{-2s}$, we also have $e^{-2s}\in V$. Hence all generators of $W$ lie in $V^\sigma$, proving the forward inclusion.

We give explicit formulas proving the reverse inclusion. Write
\[
 P=\nn{P_E^2},\quad Q=\nn{P_F^2},\quad
 X=(u_{1,7})_{1}\theta_J(u_{2,6}),
\]
and define the following weight-two vectors:
\begin{align}\label{eq:reverseaux}
 Z_0&=\tfrac14P_{1}Q-\omega-\tfrac13(P_H)_{-1}P_H-\partial P_H,\notag\\
 Y&=(Z_0)_{1}X+\tfrac43X.
\end{align}
By~\eqref{eq:Ueven}, $W=\widetilde U^\tau$. Since $\theta_H(U)=U$, $\theta_H\tau=\tau\theta_H$ and $\Phi\theta_H=\theta_J\Phi$, the subalgebra $W=\Phi(U^\tau)$ is stable under $\theta_J$. By Proposition~\ref{prop:dictionary}, $P,Q,u_{2,6},u_{1,7},\theta_Ju_{2,6}$ belong to $W$, so $X,Z_0,Y\in W$.

A direct calculation gives
\begin{equation}\label{eq:PBwitness}
 (P_{1}X)_{1}Y=-\frac{32}{3}\nn{P_EB_E}.
\end{equation}
It is an exact equality in the fixed cocycle convention. It can be evaluated solely from the weight-two generator table and the lattice vertex-operator formula. To specify oscillator handling in this calculation, use
\[
 Y(e^x,z)\alpha_i(-n)=\bigl(\alpha_i(-n)-(x,\alpha_i)z^{-n}\bigr)Y(e^x,z)
\]
and $\alpha_i(z)=Y(\alpha_i(-1)\one,z)$ together with $Y(\alpha_i(-n)a,z)=:\partial^{n-1}\alpha_i(z)Y(a,z):/(n-1)!$. Appendix~\ref{app:nested} gives explicit intermediate expansions and a three-entry Gram calculation proving that the entire difference is zero.

The remaining mixed products follow by applying $\theta_J$ and using the affine identities
\begin{align}\label{eq:otherPB}
 Q_{1}\nn{P_EB_E}&=4\nn{P_FB_E},\notag\\
 \theta_J\nn{P_EB_E}&=\nn{P_FB_F},\qquad
 P_{1}\nn{P_FB_F}=4\nn{P_EB_F}.
\end{align}
Finally the two charge exponentials are obtained by
\begin{equation}\label{eq:e2p}
 P_{-5}\bigl(P_{-3}P\bigr)=-48e^{2p},\qquad
 Q_{-5}\bigl(Q_{-3}Q\bigr)=-48e^{-2p}.
\end{equation}
For example, $P=-2\sum_{i<j}e^{\rho_i+\rho_j}$; orthogonality of the three $\rho_i$ makes~\eqref{eq:e2p} a three-term lattice calculation. All the vectors required by Proposition~\ref{prop:comparison} now belong to $W$. That proposition and the lattice-shift lemma prove equality in every weight.
\end{proof}

\section{A positive conformal-vector identity and its application}\label{sec:rationality}
In this section, we apply $\Phi^{-1}$ to the two parafermion subalgebras of $V^\sigma$ to obtain two parafermion subalgebras of $U^\tau$. We then average their conformal vectors over the compact group associated with Huang's level-$39$ affine subalgebra. The resulting positive conformal-vector identity allows us to apply the criterion of Section~\ref{sec:criterion} and prove strong rationality of $U$. The fixed point and extension results then give strong rationality of $C$ and $V$.

\subsection{The positive conformal-vector identity}
\leavevmode
\begingroup
Huang's positive conformal-vector identity uses the averaged conformal vectors of embedded copies of $K(\slc,3)$ and $K(\slc,13)$; see \cite[Section 4, (16)]{Huang}. For $k=3,13$, let $\eta_k^J\in V^\sigma$ denote the parafermion conformal vector obtained from Jiao's level-$k$ affine subalgebra by~\eqref{eq:eta}. Theorem~\ref{thm:fixed} identifies $U^\tau$ with $V^\sigma$ through the unitary isomorphism $\Phi$. We therefore set
\[
 \eta_k=\Phi^{-1}(\eta_k^J)\in U^\tau,\qquad k=3,13.
\]
Each $\eta_k$ is the conformal vector of the corresponding copy of $K(\slc,k)$ embedded in $U^\tau$ by $\Phi^{-1}$. These are the transported parafermion conformal vectors used in this section. Their central charges are $4/5$ and $8/5$, respectively. Define
\par\endgroup
\[
 q_0=\frac13\eta_3+\frac1{13}\eta_{13},\qquad
 \omega_c=\omega-\omega_L-\omega_R.
\]
\begingroup
We average over the compact real form of the right zero-mode $\slc$ action. The adjoint relations in Section~4.4 show that
\[
 \mathrm i(H_R)_0,\qquad
 (E_R)_0-(F_R)_0,\qquad
 \mathrm i\bigl((E_R)_0+(F_R)_0\bigr)
\]
are skew-adjoint operators. They span the corresponding $\mathfrak{su}_2$ action, which integrates to a compact group $G_R$ of unitary automorphisms of $U$. Define its normalized Haar average by
\[
 \mathcal E_R(v)=\int_{G_R}g(v)\,dg,
 \qquad \int_{G_R}1\,dg=1,
\]
where $dg$ is Haar measure. The integral is taken on each finite-dimensional homogeneous subspace of $U$.

Invariance of Haar measure and unitarity imply that $\mathcal E_R$ is the orthogonal projection onto the vectors fixed by $G_R$. These are precisely the vectors annihilated by $(E_R)_0$, $(F_R)_0$, and $(H_R)_0$. Consequently, the projection calculation in Section~\ref{sub:Riesen} gives
\[
 \mathcal E_R(v)=\Pi_Rv
 \qquad\bigl(v\in U_2,\ (H_R)_0v=0\bigr),
\]
with $\Pi_R$ as in~\eqref{eq:casimirprojectors}. We can therefore evaluate the averages using that explicit Casimir polynomial. Their interpretation as group averages will also express the resulting identity as a finite positive combination of conformal vectors of conjugate parafermion subalgebras.
\par\endgroup

\begin{proposition}\label{prop:positiveidentity}
The transported parafermion conformal vectors satisfy
\begin{equation}\label{eq:span}
 \mathcal E_R(q_0)=\frac{464}{22815}\omega_R+\frac2{25}\omega_c.
\end{equation}
Consequently
\begin{equation}\label{eq:positive}
 \omega=\omega_L+\frac{3403}{4563}\omega_R
       +\frac{25}{6}\mathcal E_R(\eta_3)
       +\frac{25}{26}\mathcal E_R(\eta_{13}).
\end{equation}
All coefficients in \eqref{eq:positive} are positive.
\end{proposition}
\begin{proof}
The vectors $\eta_k$ are obtained by substituting the $\slc$-triples \eqref{eq:P}--\eqref{eq:BH} into \eqref{eq:eta}, and then applying $\Phi^{-1}$ to both the momenta and the Heisenberg factors. In particular, the inverse image of a monomial is
\begin{equation}\label{eq:inverseliftfull}
 \Phi^{-1}\bigl(h_1(-n_1)\cdots h_m(-n_m)e^y\bigr)
 =\frac{(R_*^{-1}h_1)(-n_1)\cdots(R_*^{-1}h_m)(-n_m)e^{R_*^{-1}y}}
 {\gamma_*(y)\chi_*(y)}.
\end{equation}
Both parafermion conformal vectors have zero $H_R(0)$-weight. Applying the eight factors of $\Pi_R$ and collecting coefficients gives the more precise identities
\begin{align}\label{eq:separateaverages}
 \mathcal E_R(\eta_3)&=\frac2{45}\omega_R+\frac{39}{175}e+\frac9{100}f,\notag\\
 \mathcal E_R(\eta_{13})&=\frac{14}{195}\omega_R+\frac{13}{175}e+\frac{13}{20}f.
\end{align}
Here $e,f$ are the vectors recovered in \eqref{eq:recoveredvirasoro}. Appendix~\ref{app:invariants} proves these full vector identities using explicit Gram matrices, Casimir moments and a zero-norm calculation. Multiplying the first identity by $1/3$ and the second by $1/13$ makes the coefficients of $e$ and $f$ both equal to $2/25$, and their sum is $\omega_c$. This proves \eqref{eq:span}. All scalar entries needed for that verification are included in the appendix; the argument does not assume in advance that the averages belong to the displayed span.

The scalar coefficients can also be checked from
\[
 (\mathcal E_R(q_0),\omega_R)=\frac{232}{7995},\qquad
 (\mathcal E_R(q_0),\omega_c)=\frac{34}{205},
\]
since the squared norms of $\omega_R$ and $\omega_c$ are $117/82$ and $85/41$. Moreover
\[
 (q_0,\omega)
 =\frac13\frac25+\frac1{13}\frac45=\frac{38}{195},
\]
which agrees with the sum of the preceding two pairings. Finally, solve \eqref{eq:span} for $\omega_c$ and use $\omega=\omega_L+\omega_R+\omega_c$. This gives \eqref{eq:positive}.
\end{proof}

\begingroup
We now explain how the averages in~\eqref{eq:positive} give the finite positive sum required by Theorem~\ref{thm:positivecriterion}. For $k=3,13$ and $g\in G_R$, the vector $g(\eta_k)$ is the conformal vector of the conjugate of the corresponding parafermion subalgebra. This conjugate is again regular and unitary. The orbit $\{g(\eta_k)\mid g\in G_R\}$ is compact in $U_2$, viewed as a finite-dimensional real vector space. By Carath\'eodory's theorem, its convex hull is compact and hence closed; see \cite[Propositions 1.2.1--1.2.2]{Bertsekas}. Since Haar measure is normalized, $\mathcal E_R(\eta_k)$ is a limit of finite convex combinations of points in the orbit, and therefore belongs to this convex hull. Hence there exist a positive integer $N_k$, elements $g_{k,j}\in G_R$, and positive real numbers $t_{k,j}$ such that
\[
 \mathcal E_R(\eta_k)
   =\sum_{j=1}^{N_k}t_{k,j}g_{k,j}(\eta_k),
 \qquad \sum_{j=1}^{N_k}t_{k,j}=1.
\]
Zero coefficients have been omitted. The average itself need not be a Virasoro vector; the conformal vectors used in the criterion are the individual $g_{k,j}(\eta_k)$. Substituting these finite sums into~\eqref{eq:positive} gives exactly the form~\eqref{eq:positivefinite}. The relevant regular unitary vertex subalgebras are the two affine algebras and the finitely many conjugate parafermion subalgebras appearing in these sums.
\par\endgroup

\subsection{Application to the three vertex operator algebras}

\begin{corollary}\label{cor:rationality}
The vertex operator algebras $U,V$ and $C$ are strongly rational.
\end{corollary}
\begin{proof}
Corollary~\ref{cor:unitary} gives simplicity and unitarity of $U$ and $V$. The affine algebras at positive integral level and the parafermion algebras $K(\slc,3),K(\slc,13)$ are regular, and their irreducible modules have nonnegative conformal weights, with weight zero only for the vacuum. Their embeddings preserve the unitary structures, grading and translation. Their conformal vectors are annihilated by the ambient $L(1)$, as follows directly from the affine Sugawara and Heisenberg formulas. Their zero modes $(e_i)_1$ are consequently self-adjoint for the same Hermitian form on $U$. The same properties hold for their unitary conjugates. Proposition~\ref{prop:positiveidentity} and Theorem~\ref{thm:positivecriterion} therefore prove strong rationality of $U$.

The fixed point theorem of Carnahan--Miyamoto \cite{CM} gives strong rationality of $U^\tau\cong C$. Since $C=V^\sigma$ and $V$ is simple and of CFT type, the converse extension theorem for finite-group fixed point subalgebras gives strong rationality of $V$; see \cite[Example 4.2]{McRaeExt}. In this application the categorical dimension of the extension is $2$, so the nonzero-index condition is satisfied.
\end{proof}

\section{The module categories}\label{sec:categories}
\begingroup
In this section, we use the common fixed point subalgebra to compare the module categories of $U$ and $V$. We first determine the fusion products and twists of the two odd $C$-modules. We then prove that they generate a nondegenerate subcategory of $\Rep(C)$ and apply the general results on extensions to obtain a ribbon equivalence between $\Rep(U)$ and $\Rep(V)$. Finally, we apply this equivalence to the Haagerup-center identification and distinguish the two vertex operator algebras by their affine levels.

By Corollary~\ref{cor:rationality}, the vertex operator algebras $C=W$, $\widetilde U$ and $V$ are strongly rational, and $\widetilde U\cong U$. Thus their categories of ordinary modules are modular tensor categories. Let $M_s$ and $M_p$ be the $-1$ eigenspaces of $\tau$ on $\widetilde U$ and of $\sigma$ on $V$, respectively. Then
\begin{equation}
 \widetilde U=C\oplus M_s,\qquad V=C\oplus M_p,
 \qquad \mathsf a=M_s,\quad \mathsf b=M_p.
\end{equation}
All fusion products in this section are taken over $C$; we write $XY$ for $X\boxtimes_C Y$ and denote the tensor unit $C$ by $\one$ in categorical formulas. For a simple $C$-module $X$ of lowest conformal weight $h_X$, its twist is $\theta_X=e^{2\pi i h_X}$. The lattice decompositions in Section~\ref{sec:neighbor} give
\[
 M_s\subset V_{s+L_0},\qquad M_p\subset V_{p+L_0}.
\]
The subscripts record these ambient lattice cosets. All lattice modules below are regarded as $C$-modules by restriction along the conformal inclusion $C\subset V_{L_0}$.

\begin{lemma}\label{lem:twists}
The modules $\mathsf a$ and $\mathsf b$ are distinct irreducible simple currents of order two. Their twists satisfy
\[
 \theta_{\mathsf a}=\theta_{\mathsf b}=1,\qquad
 \theta_{\mathsf a\mathsf b}=-1.
\]
Consequently, the simple objects generated by $\mathsf a$ and $\mathsf b$ are precisely $\one,\mathsf a,\mathsf b,\mathsf a\mathsf b$, with fusion group $(\Z/2)^2$.
\end{lemma}
\begin{proof}
The tensor-category form of orbifold Schur--Weyl duality identifies the subcategory generated by the eigenspaces of either involution with $\Rep(\Z/2)$, the category of finite-dimensional complex representations of $\Z/2$; see \cite[Corollary 4.8]{McRaeOrb}. In particular, $M_s$ and $M_p$ are irreducible, neither is isomorphic to $C$, and
\[
 \mathsf a^2\cong\one,\qquad \mathsf b^2\cong\one.
\]
Their conformal weights are integral because they are subspaces of the vertex operator algebras $\widetilde U$ and $V$. Hence $\theta_{\mathsf a}=\theta_{\mathsf b}=1$.

To determine the twist of their fusion product, we use the lattice modules containing them. Let $L_0^*=\{x\in L_0\otimes_{\Z}\Q:(x,L_0)\subset\Z\}$ be the dual lattice. The integral overlattices $L_H$ and $L_J$ show that $s,p\in L_0^*$. The usual lattice construction gives an intertwining operator $\mathcal Y$ of type
\[
 \binom{V_{s+p+L_0}}{V_{s+L_0}\quad V_{p+L_0}}.
\]
Its restriction to $M_s$ and $M_p$ is nonzero. Indeed, the transported generator $E_L$ is a nonzero scalar multiple of $e^s$, so $e^s\in M_s$, while \eqref{eq:cubic} gives $e^p=-\frac16\nn{P_E^3}\in M_p$. On these vectors the lattice formula is
\[
 \mathcal Y(e^s,z)e^p
 =\nu z^{(s,p)}
   \exp\!\left(\sum_{n>0}\frac{s(-n)}{n}z^n\right)e^{s+p},
 \qquad \nu\in\C^\times.
\]
Here $\nu$ is a nonzero normalization constant for the intertwining operator. In particular, the coefficient of $z^{3/2}$ is $\nu e^{s+p}$, since $(s,p)=3/2$ by \eqref{eq:norms}.

The target $V_{s+p+L_0}$ is an ordinary $C$-module. By the defining property of the fusion product, the restricted intertwining operator gives a nonzero $C$-module homomorphism
\[
 \mathsf a\mathsf b\longrightarrow V_{s+p+L_0}.
\]
The source is simple because it is invertible, so this homomorphism is injective. For every $x\in L_0$,
\[
 \frac{(s+p+x)^2}{2}
 =\frac{(s+p)^2}{2}+(s+p,x)+\frac{x^2}{2}
 \equiv\frac{(s+p)^2}{2}\pmod{\Z}.
\]
The Heisenberg operators add integral weights. Thus every weight of $V_{s+p+L_0}$, and hence every weight of $\mathsf a\mathsf b$, is congruent modulo $\Z$ to
\[
 \frac{(s+p)^2}{2}
 =\frac{2+6+2(3/2)}{2}=\frac{11}{2}.
\]
It follows that $\theta_{\mathsf a\mathsf b}=-1$. In particular, $\mathsf a\mathsf b\not\cong\one$, so $\mathsf a\not\cong\mathsf b$. Since fusion is commutative up to the braiding and both objects have order two, the four displayed objects are distinct and exhaust the simple objects they generate.
\end{proof}

Let $\Tc$ be the full fusion subcategory of $\Rep(C)$ generated by $\mathsf a$ and $\mathsf b$. It is pointed, meaning that all its simple objects are invertible. Write $c_{X,Y}:X\boxtimes_C Y\longrightarrow Y\boxtimes_C X$ for the braiding. We denote by $\D$ the M\"uger centralizer of $\Tc$ in $\Rep(C)$: its objects are those $X$ for which
\[
 c_{Y,X}\circ c_{X,Y}=\mathrm{id}_{X\boxtimes_C Y}
 \quad\text{for every }Y\in\Tc.
\]
The next proposition derives the comparison of module categories directly from Lemma~\ref{lem:twists} and the algebra structures supplied by the two extensions.

\begin{proposition}\label{prop:categories}
The subcategory $\Tc$ is nondegenerate, and there are ribbon equivalences
\[
 \Rep(C)\simeq_{\rm rib}\Tc\boxtimes\D,
 \qquad \Rep(U)\simeq_{\rm rib}\D\simeq_{\rm rib}\Rep(V).
\]
\end{proposition}
\begin{proof}
By Lemma~\ref{lem:twists}, the four simple objects of $\Tc$ form the group $(\Z/2)^2$ and have twists
\[
\begin{array}{c|rrrr}
 x&\one&\mathsf a&\mathsf b&\mathsf a\mathsf b\\\hline
 \theta_x&1&1&1&-1.
\end{array}
\]
The balancing identity gives $c_{y,x}\circ c_{x,y}=\beta(x,y)\mathrm{id}_{xy}$, where $\beta(x,y)=\theta_{xy}/(\theta_x\theta_y)$. For $i,j,i',j'\in\{0,1\}$, it follows that
\[
 \beta(\mathsf a^i\mathsf b^j,\mathsf a^{i'}\mathsf b^{j'})
 =(-1)^{ij'+ji'}.
\]
This pairing has trivial radical: an object $\mathsf a^i\mathsf b^j$ centralizing both $\mathsf a$ and $\mathsf b$ must have $i=j=0$. Thus $\Tc$ is nondegenerate, and M\"uger's decomposition theorem \cite[Theorem 4.2]{Muger} gives the ribbon equivalence $\Rep(C)\simeq_{\rm rib}\Tc\boxtimes\D$.

The two extensions $\widetilde U=C\oplus\mathsf a$ and $V=C\oplus\mathsf b$ determine commutative algebra objects
\[
 \mathcal A_s=\one\oplus\mathsf a,
 \qquad \mathcal A_p=\one\oplus\mathsf b
\]
in $\Tc$ by \cite[Theorem 3.2]{HKL}. The orbifold equivalences used in the proof of Lemma~\ref{lem:twists} identify the subcategory generated by either odd summand with $\Rep(\Z/2)$. The product of that summand with itself is nonzero by simplicity of the extension. After rescaling, each algebra is therefore the algebra of functions on $\Z/2$ with its translation action. Hence both algebras are connected, meaning $\operatorname{Hom}(\one,\mathcal A_i)\cong\C$ for $i=s,p$, and separable; see \cite[Example 3.3(i)]{DMNO}.

Write $\FPdim$ for Frobenius--Perron dimension; the dimension of a fusion category is the sum of the squared dimensions of its simple objects. For a commutative algebra $\mathcal A$ in a braided fusion category $\mathcal E$, let $\mathcal E^0_{\mathcal A}$ denote the category of local right $\mathcal A$-modules, whose action is invariant under double braiding with $\mathcal A$. Since $\FPdim(\Tc)=4$ and $\FPdim(\mathcal A_s)=\FPdim(\mathcal A_p)=2$, the dimension formula \cite[Corollary 3.32]{DMNO} gives
\[
 \FPdim(\Tc^0_{\mathcal A_i})
 =\frac{\FPdim(\Tc)}{\FPdim(\mathcal A_i)^2}
 =\frac4{2^2}=1\qquad(i=s,p).
\]
Thus $\Tc^0_{\mathcal A_s}\simeq\mathrm{Vec}\simeq\Tc^0_{\mathcal A_p}$. Both extension algebras lie in the $\Tc$ factor of $\Rep(C)$, so
\[
 \Rep(C)^0_{\mathcal A_i}
 \simeq\Tc^0_{\mathcal A_i}\boxtimes\D
 \simeq\D\qquad(i=s,p).
\]
By \cite[Theorem 3.4]{HKL} and \cite[Theorem 3.65]{CKMExt}, these are braided equivalences identifying $\Rep(\widetilde U)$ and $\Rep(V)$ with $\D$. The extension algebras have trivial twist. The local-module twists are inherited from $\Rep(C)$ and agree with $e^{2\pi iL(0)}$ for the extensions, since all inclusions preserve the conformal vector. The equivalences are therefore ribbon equivalences. Finally, $\widetilde U\cong U$ gives the assertion for $\Rep(U)$.
\end{proof}

The classification of pointed braided fusion categories, applied to the fusion group and twists in Lemma~\ref{lem:twists}, gives
\[
 \Tc\simeq_{\rm rib}Z(\mathrm{Vec}_{\Z/2});
\]
see \cite[Theorem 8.4.12]{EGNO}. Here the ribbon structures are the canonical positive ones: all four simple objects of $\Tc$ have categorical dimension one, by the orbifold equivalences used in Lemma~\ref{lem:twists}.

The following corollary transfers the Haagerup-center identification from either realization to the other. The rank of a fusion category is the number of isomorphism classes of its simple objects.

\begin{corollary}
If either $\Rep(U)$ or $\Rep(V)$ is ribbon equivalent to $Z(\mathcal H_3)$, then both are, and
\[
 \Rep(C)\simeq_{\rm rib}
 Z(\mathrm{Vec}_{\Z/2})\boxtimes Z(\mathcal H_3).
\]
In particular, $\Rep(C)$ has rank $48$ when $Z(\mathcal H_3)$ has rank $12$.
\end{corollary}
\begin{proof}
Proposition~\ref{prop:categories} identifies both module categories with $\D$. The preceding identification of $\Tc$ then gives the displayed formula. Simple objects of a Deligne tensor product are indexed by pairs of simple objects of its factors. Since $\Tc$ has four simple objects, the rank is $4\cdot12=48$.
\end{proof}

This comparison is stronger than equality in the Witt group. Since $\Tc$ is a Drinfeld center, its Witt class is zero \cite[Section 5.1]{DMNO}, and the factorization in Proposition~\ref{prop:categories} gives
\[
 [\Rep(C)]=[\Rep(U)]=[\Rep(V)].
\]
Here brackets denote the class in the Witt group of nondegenerate braided fusion categories. If $\D\simeq Z(\mathcal H_3)$, all these classes are zero.

Finally, a vertex operator algebra isomorphism preserves both the Lie algebra structure and the invariant bilinear form on the weight-one space, and hence the levels of its simple factors. For $U$ and $V$, these levels are $1,39$ and $3,13$, respectively, by the full weight-one calculations in \cite[Section 3, (15)]{Huang} and \cite[Proposition 4.6]{Jiao}. We therefore have the following proposition.

\begin{proposition}\label{prop:distinct}
The vertex operator algebras $U$ and $V$ are not isomorphic as conformal vertex operator algebras.
\end{proposition}

The full weight-one equalities used in Proposition~\ref{prop:distinct} are inputs from the original constructions. They do not follow solely from the displayed generators, since products of weight-two generators can produce weight-one elements. They are not needed for Theorem~\ref{thm:fixed} or Proposition~\ref{prop:categories}. The factorization of $\Rep(C)$ concerns module categories; it does not assert a tensor product decomposition of the vertex operator algebra $C$, or an automorphism of $C$ interchanging $M_s$ and $M_p$.
\par\endgroup

\section{The associated conformal nets}\label{sec:nets}
We prove that the conformal nets associated with $U$, $V$ and $C$ are completely rational. The proof compares the constructions using bounded localized and smeared vertex operators and applies a finite-index theorem for conformal subnets. We then show that the common fixed point isomorphism of vertex operator algebras induces the corresponding isomorphism of fixed point nets.

In this section $W$ denotes any one of $U,V,C$. We distinguish the net $\mathcal A_W^{\rm loc}$ constructed from bounded localized vertex operators from the net $\mathcal A_W^{\rm sm}$ constructed from smeared vertex operators. We will show that they are isomorphic for the vertex operator algebras considered here.

Put $E=V_{E_8}$. The positive lattice characters $\kappa_H$ and $\kappa_J$ extend from $L=E_6\oplus A_2$ to its index-three overlattice $E_8$. Explicitly, for $x\in E_8$ the extension is
\[
 \widetilde\kappa(x)=\kappa(3x)^{1/3},
\]
where the positive cube root is taken. Transporting the standard lattice unitary structure by this character realizes each of $U,V$, and $C\subset V$, as a unitary conformal subalgebra of a copy of $E$. The resulting unitary lattice VOAs are isomorphic. Thus the different normalizations of Section~\ref{sec:lattice} cause no change in the analytic properties used below.

\begin{theorem}\label{thm:nets}
For $W=U,V$, or $C$, the vertex operator algebra $W$ is strongly local and has bounded localized vertex operators. Both $\mathcal A_W^{\rm loc}$ and $\mathcal A_W^{\rm sm}$ are completely rational, and
\[
 \mathcal A_W^{\rm loc}\cong\mathcal A_W^{\rm sm}.
\]
The subnet of $\mathcal A_E^{\rm loc}$ generated by $W$ has finite Jones index.
\end{theorem}
\begin{proof}
The lattice VOA $E\cong L_{\mathfrak e_8}(1,0)$ is strongly local: it is generated by its weight-one space, so \cite[Theorem 8.3]{CKLW} applies. Strong locality passes to its unitary subalgebras by \cite[Theorem 7.1]{CKLW}. Bounded localized vertex operators for $E$ follow from \cite[Theorem 7.4]{TenerFields}, and pass to $W$ by \cite[Proposition 4.16]{TenerFields}.

Let $\mathcal B_W\subset\mathcal A_E^{\rm loc}$ be the subnet generated by the bounded localized vertex operators of $W$ in the Hilbert space completion of $E$. Its vacuum restriction to the closure of $W$ is $\mathcal A_W^{\rm loc}$. By Corollary~\ref{cor:rationality}, $W$ is regular. Its commutant in $E$ is trivial. Indeed, if $a\in C_E(W)$, then $\omega\in W$ gives
\[
 L(0)a=\omega_1a=0,
\]
so $a\in E_0=\C\one$. Consequently
\[
 C_E(W)=\C\one,
\]
which is regular. All the hypotheses of \cite[Corollary 4.24]{TenerFusion} are therefore satisfied, and
\begin{equation}\label{eq:netindex}
 [\mathcal A_E^{\rm loc}(I):\mathcal B_W(I)]<\infty
\end{equation}
for every nonempty proper interval $I\subset S^1$, where $S^1$ is the unit circle and the bracket denotes the Jones index of the inclusion of local von Neumann algebras.

The $E_8$ lattice net is completely rational by \cite[Corollary 3.19]{DXNets}. The standard identification with the level-one $E_8$ loop-group net, together with \cite[Example 8.7]{CKLW} and \cite[Corollary 7.7]{TenerFields}, identifies $\mathcal A_E^{\rm loc}$ with this lattice net up to isomorphism. Longo's finite-index theorem \cite[Theorem 24]{Longo} now implies complete rationality of $\mathcal A_W^{\rm loc}$.

For the ambient affine VOA $E$, the two net constructions are isomorphic by \cite[Corollary 7.7]{TenerFields}. The subalgebra comparison in \cite[Lemma 7.6]{TenerFields} then gives $\mathcal A_W^{\rm loc}\cong\mathcal A_W^{\rm sm}$. Hence the smeared-field net is completely rational as well. This comparison is an isomorphism of conformal nets; no equality of the two concrete families of operator algebras is required.
\end{proof}

From now on write $\mathcal A_W=\mathcal A_W^{\rm sm}$. The common fixed point result has the following net counterpart.

\begin{corollary}\label{cor:netfixed}
There are isomorphisms of completely rational conformal nets
\[
 \mathcal A_U^{\langle\tau\rangle}\cong\mathcal A_C
       \cong\mathcal A_V^{\langle\sigma\rangle}.
\]
\end{corollary}
\begin{proof}
The lattice-character involutions $\tau$ and $\sigma$ preserve the positive forms, and hence induce unitary automorphisms of the associated nets. The fixed point compatibility theorem \cite[Proposition 7.6]{CKLW} gives
\[
 \mathcal A_{U^\tau}=\mathcal A_U^{\langle\tau\rangle},\qquad
 \mathcal A_{V^\sigma}=\mathcal A_V^{\langle\sigma\rangle}.
\]
The isomorphism $\Phi:U^\tau\longrightarrow C=V^\sigma$ preserves the unitary structures by Proposition~\ref{prop:dictionary}. Its extension to the Hilbert space completions intertwines the smeared vertex operators and therefore the associated nets. Theorem~\ref{thm:nets} gives complete rationality.
\end{proof}

\begin{remark}
Complete rationality implies that the Doplicher--Haag--Roberts (DHR) sectors form a unitary modular tensor category \cite{KLM}. The results of this section do not require an equivalence between that category and the VOA module category. In particular, complete rationality alone does not supply the identification of the DHR category with the Haagerup center. The VOA module-category comparison in Section~\ref{sec:categories} is independent of this analytic identification.
\end{remark}

\section{Further questions}\label{sec:remaining}
We discuss the scope of the comparison and the additional input needed to identify the module categories with the Haagerup center. We also consider possible extensions of the common fixed point construction and of the positive conformal-vector method to other rational vertex operator algebras.

The preceding results establish the equality of the Riesen and Huang constructions after an explicit lattice isomorphism, and identify their common order-two fixed point subalgebra with that of Jiao's construction. The calculation of the transported parafermion conformal vectors also gives the rationality argument directly in this realization. None of these steps requires a classification of the irreducible modules.

The identification of the resulting module category with $Z(\mathcal H_3)$ uses the additional branching and categorical results in the original constructions. The present comparison transfers such an identification between the algebras. It does not by itself supply the branching calculation or classify all irreducible modules. The full weight-one spaces and common character quoted in the Introduction are inputs from the original constructions; Proposition~\ref{prop:distinct} uses the weight-one results. These inputs have not been rederived from the finite comparison identities in this paper.

It would be interesting to determine whether similar fixed point comparisons occur for other rational vertex operator algebras with the same character. One may also ask when the parafermion subalgebras in a common fixed point subalgebra yield a positive conformal-vector identity to which Theorem~\ref{thm:positivecriterion} applies.

\appendix
\section{Finite lattice calculations}\label{app:checks}
This appendix supplies the coefficient formulas and finite calculations used in the comparison. We give the lattice mode recursions and Hermitian pairings, the intermediate expressions for the common fixed point identity, and the complete coefficient tables and Gram matrices for the Riesen--Huang comparison and the parafermion averages. These data prove the required vector identities by computing the norm of the entire difference and permit the calculations to be checked independently of the accompanying programs.

\subsection{Coefficient formulas and Hermitian pairings}\label{app:kernel}
Let $M$ be one of the positive definite even lattices considered in the paper, with ordered basis $\alpha_i$, Gram matrix $G$ and the chosen cocycle $\epsilon$. Write $x=\sum_i x_i\alpha_i$ for lattice coordinates. Introduce commuting symbols $a_{i,k}$ for $\alpha_i(-k)$, $k\ge1$. A Fock monomial is $Pe^x$, where $x\in M$ and $P$ is a monomial in these symbols; its weight is $x^2/2+\sum k\deg_{a_{i,k}}P$. Put
\[
 E_x(t)=\exp\!\left(\sum_{i,k\ge1}\frac{x_i}{k}a_{i,k}t^k\right),\qquad
 E_x[d]=\sum_{\substack{m_{i,k}\ge0\\\sum k m_{i,k}=d}}
 \prod_{i,k}\frac{x_i^{m_{i,k}}a_{i,k}^{m_{i,k}}}{k^{m_{i,k}}m_{i,k}!}.
\]
Only finitely many $m_{i,k}$ are nonzero, and $E_x[d]=0$ for $d<0$. Here $[t^d]$ denotes the coefficient of $t^d$. Commuting the target oscillators through a lattice exponential gives
\begin{equation}\label{eq:extraction}
 (e^x)_n(Pe^y)=\epsilon(x,y)e^{x+y}
 [t^{-n-1-(x,y)}]E_x(t)
 P\bigl(a_{i,k}-(x,\alpha_i)t^{-k}\bigr).
\end{equation}
Oscillators on the first vector are handled recursively by
\begin{align}\label{eq:oscillatorrecursion}
 (\alpha_i(-k)a)_n b={}&
 \sum_{j<0}(-1)^{k-1}\binom{j+k-1}{k-1}
       \alpha_i(j)a_{n-j-k}b\notag\\
 &+\sum_{j\ge0}(-1)^{k-1}\binom{j+k-1}{k-1}
       a_{n-j-k}\alpha_i(j)b.
\end{align}
The binomial coefficients have their polynomial meaning for negative upper arguments. The Heisenberg actions needed here are
\[
 \alpha_i(0)(Pe^x)=(\alpha_i,x)Pe^x,\qquad
 \alpha_i(j)(Pe^x)=j\sum_lG_{il}\frac{\partial P}{\partial a_{l,j}}e^x
 \quad(j>0).
\]
For $k>0$, we have $\alpha_i(-k)(Pe^x)=a_{i,k}Pe^x$. In a calculation with output weight $N\ge0$, the first sum in \eqref{eq:oscillatorrecursion} may be restricted to $-N\le j<0$. In the second sum, $j$ does not exceed the largest oscillator index of the target, apart from $j=0$. These bounds and the explicit formula for $E_x[d]$ give a finite calculation for every entry below.

The Hermitian form is linear in the first variable. If
$P=\prod_{a=1}^r a_{i_a,n_a}$ and $Q=\prod_{b=1}^s a_{j_b,m_b}$, its value is zero unless $x=y$ and $r=s$; otherwise
\begin{equation}\label{eq:fockpairing}
 (Pe^x,Qe^x)_\kappa=
 \kappa(x)\sum_{\pi\in S_r}\prod_{a=1}^r
 n_a\,\delta_{n_a,m_{\pi(a)}}G_{i_a,j_{\pi(a)}}.
\end{equation}
Here $S_r$ is the permutation group on $r$ letters and $\delta$ is the Kronecker delta. Coefficients in the second vector are conjugated. The sum for $r=0$ is $1$. The positive characters $\kappa_H,\kappa_J$ are specified in Section~\ref{sec:lattice}.
Formula~\eqref{eq:fockpairing} follows by moving each positive Heisenberg mode to the right. It is positive definite because $G$ is positive definite and $\kappa(x)>0$.

All entries are reduced using $r^2=13$, $\xi^2=-3$ and $\overline\xi=-\xi$. Thus the comparisons are exact algebraic identities. In particular, if $A,B$ and a scalar $c$ are specified, the three pairings $(A,A)$, $(A,B)$ and $(B,B)$ prove $A=cB$ whenever
\begin{equation}\label{eq:normtest}
 (A,A)-\overline c(A,B)-c(B,A)+|c|^2(B,B)=0.
\end{equation}

\subsection{The nested product for the common fixed point}\label{app:nested}
Use Jiao's coordinates and Hermitian form. Let $P,Q,X,Z_0,Y$ be the vectors in the proof of Theorem~\ref{thm:fixed}, and put $R_0=P_1X$. We give shorter intermediate expressions to make the calculation of $R_{0,1}Y$ explicit. In this subsection $R_{0,1}$ denotes the mode of the vector $R_0$, not a lattice matrix.

Put $I_0=\{1,2,4,5,7,8\}$ and
\[
 d=\left(2-\frac{2r\xi}{9}\right)(\alpha_1+\alpha_4)
 +\left(2+\frac{2r\xi}{9}\right)(\alpha_2+\alpha_5)
 -\left(4+\frac{4\xi}{9}\right)\alpha_7
 -\left(4-\frac{4\xi}{9}\right)\alpha_8.
\]
Then \eqref{eq:extraction}--\eqref{eq:oscillatorrecursion} give
\begin{equation}\label{eq:explicitP1X}
 R_0=\sum_{i\in I_0}r_i e^{\alpha_i-\alpha_6}+d(-1)e^{p+\alpha_6},
\end{equation}
where
\[
\begin{array}{c|cccc}
 i&1,4&2,5&7&8\\\hline
 r_i&-2-2r\xi/3&-1+r\xi/3&4-4\xi/3&2+2\xi/3.
\end{array}
\]
Let $\rho_1,\rho_2,\rho_3$ be the three orthogonal roots defining $P_E$ in Section~\ref{sec:lattice}, in the order used there, and set
\[
 \mu_1=\rho_3-p-\alpha_6,\quad
 \mu_2=\rho_2-p-\alpha_6,\quad
 \mu_3=\rho_1-p-\alpha_6,
 \quad\mu_4=\mu_1+\alpha_2,\quad\mu_5=\mu_2+\alpha_1.
\]
Write $k_0=-1-r\xi/3$ and $k_1=-1/2+r\xi/6$. The second intermediate vector is
\begin{align}\label{eq:explicitY}
 Y={}&-\frac12d(-1)\sum_{j=1}^3e^{\rho_j+\alpha_6}
 +\sum_{j=1}^3\left[\left(-2+\frac{2\xi}{3}\right)e^{\mu_j+\alpha_7}
                  +\left(-1-\frac\xi3\right)e^{\mu_j+\alpha_8}\right]\notag\\
 &+k_0\bigl(e^{\mu_1+\alpha_1}+e^{\mu_1+\alpha_4}\bigr)
  +k_1\bigl(e^{\mu_3+\alpha_2}+e^{\mu_3+\alpha_5}\bigr)\notag\\
 &+k_1(\alpha_2-\alpha_5)(-1)e^{\mu_4}
  +k_0(\alpha_1-\alpha_4)(-1)e^{\mu_5}.
\end{align}
The first expansion has twelve Fock monomials and the second has thirty-two. They are obtained from $X=(u_{1,7})_1\theta_J(u_{2,6})$ and $Y=(Z_0)_1X+4X/3$, without using the claimed final product. The following intermediate norms provide additional checks on their normalization:
\[
\begin{array}{c|rrrrrr}
 A&P&Q&X&Z_0&Y&R_0\\\hline
 (A,A)_{\kappa_J}&12&12&988/9&14&104&416/3.
\end{array}
\]
Now set $A=(R_0)_1Y$ and $B=(P_E)_{-1}B_E$. Substitution of \eqref{eq:explicitP1X}--\eqref{eq:explicitY} into \eqref{eq:extraction}--\eqref{eq:fockpairing} gives
\begin{equation}\label{eq:nestedgram}
 \begin{pmatrix}(A,A)&(A,B)\\(B,A)&(B,B)\end{pmatrix}_{\kappa_J}
 =\begin{pmatrix}13312/3&-416\\-416&39\end{pmatrix}.
\end{equation}
In particular
\[
 \left\|A+\frac{32}{3}B\right\|^2
 =\frac{13312}{3}-\frac{64\cdot416}{3}
       +\frac{1024\cdot39}{9}=0.
\]
Positive definiteness proves \eqref{eq:PBwitness}, including all momentum and oscillator components. This argument uses neither a dimension estimate for the fixed point subalgebra nor a test on selected coordinates.

\subsection{The invariant vectors and the two averages}\label{app:invariants}
We give a finite Gram-matrix calculation for both the Riesen--Huang identification and the positive conformal-vector identity. All vectors in this subsection are in Huang's coordinates and use $(\cdot,\cdot)_{\kappa_H}$. Set
\[
 (\nu_1,\ldots,\nu_7)=(\alpha_7,\alpha_8,\alpha_1,\alpha_2,\alpha_3,\alpha_4,\alpha_5),
 \qquad \mathsf S(x)=e^x+\kappa_H(x)e^{-x}.
\]
We retain Riesen's labels $b_0,b_1,b_2,b_3$ for his four weight-two vectors; these subscripts are unrelated to the root indices. The transported vectors $B_k=\Psi(b_k)$, $0\le k\le3$, are specified by
\begin{equation}\label{eq:Btables}
 B_k=\sum_{1\le i\le j\le7}q_{ij}^{(k)}\nu_i(-1)\nu_j(-1)\one
       +\sum_{\Gamma}g_\Gamma^{(k)}\sum_{x\in\Gamma}\mathsf S(x).
\end{equation}
Tables~\ref{tab:Bquadratic} and \ref{tab:Bexponential} give all nonzero coefficient rows. An index pair or vector appearing in a grouped row has the same displayed coefficient. All omitted coefficients are zero. The exponential coefficients already include the cocycle correction and the lift $\Psi$; no further sign convention is required. These are the transported tables of \cite[Appendix C]{RiesenLong}, included here so that the comparisons do not require a separate data file.

\begin{table}[!htbp]
\centering\small\renewcommand{\arraystretch}{1.15}
\caption{Quadratic coefficients of $B_0,B_1,B_2,B_3$}\label{tab:Bquadratic}
\begin{tabular}{@{}lrrrr@{}}\toprule
$(i,j)$ & $q^{(0)}$ & $q^{(1)}$ & $q^{(2)}$ & $q^{(3)}$\\\midrule
$(1,1),\ (1,2),\ (2,2)$ & $0$ & $0$ & $\frac{35}{3}$ & $- \frac{35}{6}$\\
$(1,3),\ (2,7)$ & $0$ & $0$ & $- \frac{10}{3}$ & $\frac{5}{6}$\\
$(1,4),\ (2,6)$ & $0$ & $0$ & $- \frac{20}{3}$ & $\frac{8}{3}$\\
$(1,5),\ (2,5)$ & $0$ & $0$ & $-9$ & $\frac{9}{2}$\\
$(1,6),\ (2,4)$ & $0$ & $0$ & $- \frac{28}{3}$ & $\frac{16}{3}$\\
$(1,7),\ (2,3)$ & $0$ & $0$ & $- \frac{20}{3}$ & $\frac{25}{6}$\\
$(3,3),\ (7,7)$ & $\frac{1}{3}$ & $\frac{4}{9}$ & $\frac{5}{3}$ & $- \frac{5}{6}$\\
$(3,4),\ (4,4),\ (6,6),\ (6,7)$ & $\frac{1}{3}$ & $\frac{2}{9}$ & $\frac{8}{3}$ & $- \frac{4}{3}$\\
$(3,5),\ (5,7)$ & $0$ & $- \frac{4}{9}$ & $2$ & $-1$\\
$(3,6),\ (4,7)$ & $- \frac{1}{3}$ & $- \frac{5}{9}$ & $\frac{4}{3}$ & $- \frac{2}{3}$\\
$(3,7)$ & $- \frac{2}{3}$ & $- \frac{1}{9}$ & $\frac{2}{3}$ & $- \frac{1}{3}$\\
$(4,5),\ (5,6)$ & $0$ & $\frac{1}{9}$ & $4$ & $-2$\\
$(4,6)$ & $- \frac{2}{3}$ & $- \frac{2}{9}$ & $\frac{8}{3}$ & $- \frac{4}{3}$\\
$(5,5)$ & $0$ & $\frac{1}{3}$ & $3$ & $- \frac{3}{2}$\\
\bottomrule\end{tabular}
\end{table}
\begin{table}[!htbp]
\centering\small\renewcommand{\arraystretch}{1.15}
\caption{Exponential coefficients of $B_0,B_1,B_2,B_3$}\label{tab:Bexponential}
\begin{tabular}{@{}lrrrr@{}}\toprule
$\Gamma$ & $g^{(0)}$ & $g^{(1)}$ & $g^{(2)}$ & $g^{(3)}$\\\midrule
$\{-\nu_1+\nu_7,\ -\nu_2+\nu_3\}$ & $0$ & $0$ & $4$ & $-3$\\
$\{-\nu_2+\nu_7,\ -\nu_1+\nu_3\}$ & $0$ & $0$ & $0$ & $1$\\
$\{-\nu_1+\nu_6,\ -\nu_2+\nu_4\}$ & $0$ & $0$ & $3$ & $-2$\\
$\{-\nu_2+\nu_6,\ -\nu_1+\nu_4\}$ & $0$ & $0$ & $1$ & $0$\\
$\{-\nu_1-\nu_2+\nu_i+\nu_{i+1}:3\le i\le6\}$ & $0$ & $0$ & $-1$ & $1$\\
$\{\nu_5-\nu_7,\ \nu_4-\nu_7\}$ & $0$ & $\frac{5}{9}$ & $0$ & $0$\\
$\{-\nu_1+\nu_5,\ -\nu_2+\nu_5\}$ & $0$ & $0$ & $2$ & $-1$\\
$\{\nu_4-\nu_6\}$ & $1$ & $0$ & $0$ & $0$\\
$\{\nu_3-\nu_5\}$ & $0$ & $1$ & $0$ & $0$\\
$\{\nu_3-\nu_6\}$ & $0$ & $\frac{8}{9}$ & $0$ & $0$\\
$\{\nu_3-\nu_7\}$ & $1$ & $- \frac{1}{3}$ & $0$ & $0$\\
$\{\nu_3+\nu_4-\nu_6-\nu_7\}$ & $-1$ & $\frac{8}{9}$ & $0$ & $0$\\
\bottomrule\end{tabular}
\end{table}

For example, \eqref{eq:fockpairing} applied to these tables gives
\[
 \bigl((B_i,B_j)\bigr)_{0\le i,j\le3}=
 \begin{pmatrix}
 10&-2/3&8&-4\\
 -2/3&104/9&32/3&-16/3\\
 8&32/3&1965&-2245/2\\
 -4&-16/3&-2245/2&2805/4
 \end{pmatrix}=:G_B.
\]
Let $\mathcal B=(\omega_L,\omega_R,B_0,B_1,B_2,B_3)$, regarded as an ordered list, and let $G_{\mathcal B}$ be its Gram matrix. The mixed pairings with $B_i$ are zero, so
\begin{equation}\label{eq:invariantgram}
 G_{\mathcal B}=\operatorname{diag}(1/2,117/82,G_B),\qquad
 \det G_{\mathcal B}=9586850.
\end{equation}
The same coefficient formulas give, for $0\le i\le3$,
\[
 \sum_{a\in\{E_L,F_L,H_L,E_R,F_R,H_R\}}
 \sum_{n=0}^{2}\|a_nB_i\|^2=0.
\]
Thus the $B_i$ belong to the affine commutant; higher positive modes vanish by weight. Also
\[
 (\omega_c,\omega_c)=85/41,
 \qquad ((\omega_c,B_i))_{i=0}^3=(2,8/3,45,-45/2).
\]
Applying \eqref{eq:normtest} with $G_B$ proves
\begin{equation}\label{eq:omegaCB}
 \omega_c=(8B_0+9B_1+2B_2+2B_3)/41.
\end{equation}
These tests use the specified four-dimensional space, without asserting that it is the full ambient commutant in weight two.

Put $a_0=\Pi_L(u_1\theta_Hu)$, $a_3=\eta_3$ and $a_{13}=\eta_{13}$, with the latter two vectors transported by \eqref{eq:inverseliftfull}. Their pairings with the six vectors of $\mathcal B$ are
\begin{equation}\label{eq:sourcepairings}
\begin{array}{c|ccc}
 b&(a_0,b)&(a_3,b)&(a_{13},b)\\\hline
 \omega_L&3&0&0\\
 \omega_R&585/41&13/205&21/205\\
 B_0&26/5&13/50&143/150\\
 B_1&164/15&41/75&91/225\\
 B_2&136+8r&(34+2r)/5&(260-26r)/15\\
 B_3&-68-8r&(-17-2r)/5&(-130+26r)/15.
\end{array}
\end{equation}
For a further calculation using only the specified affine modes, set $D_0=1$ and
\[
 D_j=\left(1-\frac{\mathcal C_R}{j(j+2)}\right)D_{j-2}
 \quad(j=2,4,\ldots,16).
\]
The following are scalar moments, not norms of the intermediate vectors:
\begin{equation}\label{eq:projectionmoments}
\begin{array}{c|rrr}
 j&(a_0,D_ja_0)&(a_3,D_ja_3)&(a_{13},D_ja_{13})\\\hline
 0&672&2/5&4/5\\
 2&-1044&-152/75&-958/325\\
 4&2700&1429/225&2727/325\\
 6&-4164&-2432/225&-1836/325\\
 8&20988/5&11357/1125&2481/1625\\
 10&-7092/5&-4438/1125&-54/1625\\
 12&6516/35&524/7875&4692/11375\\
 14&6516/35&524/7875&4692/11375\\
 16&6516/35&524/7875&4692/11375.
\end{array}
\end{equation}
Every entry in \eqref{eq:sourcepairings}--\eqref{eq:projectionmoments} is obtained by the finite recursions \eqref{eq:extraction}--\eqref{eq:fockpairing}. The initial vector $u$ has three terms; the parafermion inputs are the two explicit expressions \eqref{eq:eta} transported by \eqref{eq:inverseliftfull}. Thus these tables specify both the calculation and all scalar values needed in the argument.

The operator $D_{16}=\Pi_R$ is an orthogonal projection on the zero-Cartan-weight space under consideration. In particular the last line of \eqref{eq:projectionmoments} gives the squared norms of $\Pi_Ra_0,\Pi_Ra_3,\Pi_Ra_{13}$. Write
\[
 e_*=\frac{2B_0+6B_1}{31}+\frac{2-r}{93}B_2+\frac{2-2r}{93}B_3,
 \qquad f_*=\omega_c-e_*.
\]
Define three candidate vectors independently of the projected inputs:
\begin{align*}
 T_0&=6\omega_L+10\omega_R+\tfrac{156}{35}e_*+\tfrac95f_*,\\
 T_3&=\tfrac2{45}\omega_R+\tfrac{39}{175}e_*+\tfrac9{100}f_*,\\
 T_{13}&=\tfrac{14}{195}\omega_R+\tfrac{13}{175}e_*+\tfrac{13}{20}f_*.
\end{align*}
Each $T_i$ is annihilated by the right affine zero modes. If $c_i$ is its coordinate column in $\mathcal B$ and $s_i$ the corresponding column of \eqref{eq:sourcepairings}, the displayed matrices give
\[
\begin{array}{c|ccc}
 i&0&3&13\\\hline
 c_i^TG_{\mathcal B}c_i=c_i^Ts_i&(6516/35)&(524/7875)&(4692/11375).
\end{array}
\]
All coefficients in these columns are real. Since $(\Pi_Ra_i,T_i)=(a_i,T_i)$, equations \eqref{eq:invariantgram}--\eqref{eq:projectionmoments} imply
\[
 \|\Pi_Ra_i-T_i\|^2=(a_i,\Pi_Ra_i)-2c_i^Ts_i+c_i^TG_{\mathcal B}c_i=0.
\]
Consequently $\widetilde u=T_0$ and $\mathcal E_R(\eta_k)=T_k$. The first equality, together with \eqref{eq:recoveredvirasoro}, gives $e=e_*$ and $f=f_*$. The other two are precisely \eqref{eq:separateaverages}. This proves the full vector identities; no assumption that the projected inputs lie in the span of $\mathcal B$ is used.

For completeness, the additional Virasoro calculations in Section~\ref{sec:lattice} can also be reduced to small scalar tests. The tables give $\|e\|^2=35/31$, $\|f\|^2=1200/1271$ and $(e,f)=0$. Formula~\eqref{eq:extraction} gives the following triples $(\|A\|^2,(A,B),\|B\|^2)$:
\[
\begin{array}{c|c|ccc}
 A&B&\|A\|^2&(A,B)&\|B\|^2\\\hline
 e_1e&e&140/31&70/31&35/31\\
 f_1f&f&4800/1271&2400/1271&1200/1271\\
 e_1u&u&2028/961&156/31&12\\
 f_1u&u&388800/1615441&2160/1271&12.
\end{array}
\]
The zero-norm test proves the stated products. The vectors $e,f$ are quasiprimary, as is visible from \eqref{eq:Btables}, which has only quadratic oscillators and norm-four exponentials. Skew-symmetry gives $e_2e=f_2f=0$ and $e_0e=\partial e$, $f_0f=\partial f$. Their third products are their squared norms. Finally $e+f=\omega_c$ and their membership in the affine commutant give all mixed singular products equal to zero.

\subsection{Independent verification and external inputs}\label{app:programs}
The formulas and tables in this appendix are part of the proof. The programs provide an independent way to repeat the finite arithmetic. The files \nolinkurl{verify_comparison.py}, \nolinkurl{lattice_mode_engine.py} and \nolinkurl{verify_repair.py} verify the order-three lattice data, cocycles, affine brackets, relations involving $\theta_H$ and $\theta_J$, and the common fixed point identities. The files \nolinkurl{riesen_coset_data.py} and \nolinkurl{verify_additional.py} verify the invariant-vector comparisons, Virasoro products and averages. Their full expansions are recorded in \nolinkurl{additional_coefficients.json}. The additional scripts in \nolinkurl{paper_checks/} reproduce the intermediate expressions, Gram matrices and scalar tables printed here. All calculations use exact algebraic coefficients and assertions enabled. The exported coefficient tables use the root indices $1,\ldots,8$ of the paper; internally, Python array position $i-1$ represents $\alpha_i$. The auxiliary basis $\nu_1,\ldots,\nu_7$ is indexed in the same way.

The general rationality theorem, the all-weight lattice-shift argument and the tensor-category arguments are proved in the text. The full characters, full weight-one spaces and Haagerup-center identifications remain cited results from the original constructions. The appendix verifies the new finite calculations used in this paper; it does not replace those cited theorems by a low-weight numerical test.

\section*{Acknowledgements and AI disclosure}
The authors thank the researchers whose contemporaneous constructions made the comparison in this paper possible. The manuscript grew out of discussions aimed at simplifying several AI-assisted proofs of the Haagerup conformal field theory. AI systems were used for literature searches, proposing and checking proof strategies, exact lattice calculations, and manuscript preparation. The authors are responsible for the mathematical content.

\end{document}